\documentclass[a4paper,11pt]{article}
\usepackage{amssymb}
\usepackage{amsmath}
\usepackage{amscd}
\usepackage[francais]{babel}
\usepackage[all]{xypic}
\usepackage[all]{xy}
\usepackage{url}

\usepackage{comment}

\newtheorem{theo}{Theor\`eme}[section]
\newtheorem{prop}[theo]{Proposition}
\newtheorem{lem}[theo]{Lemme}
\newtheorem{cor}[theo]{Corollaire}
\newtheorem{defi}[theo]{D\'efinition}

\renewcommand{\thetheo} {\arabic{section}.\arabic{theo}}

\def \kbar {{\bar k}}

\def \dem {\paragraph{D\'emonstration : }}

\def \rem {\paragraph{Remarque : }}
\def \rems {\paragraph{Remarques : }}

\def \hyp{{\bf H}}
\def \Romannumeral #1 {\expandafter\uppercase\expandafter {\romannumeral #1} }
\def \br {{\rm{Br\,}}}
\def\Brr{\mathrm{Br}}
\def\brr{\mathrm{Br}}
\def \sc{{\rm sc}}
\def \ses{{\rm ss}}

\def \bra {{\rm{Br_1\,}}}
\def \brnr {{\rm{Br_{\rm nr}\,}}}

\def \pic {{\rm {Pic\,}}}
\def \div {{\rm{Div\,}}}
\def \gal {{\rm{Gal}}}

\def \calo {{\cal O}}
\def \spec {{\rm{Spec\,}}}

\def \Hom {{\rm {Hom}}}
\def \fhom {{\mathcal Hom}}
\def \ext {{\rm {Ext}}}

\def \into{{\tt {Inn}}\,}

\def \Z {{\bf Z}}
\def \Q {{\bf Q}}
\def \F {{\bf F}}

\def \NN {{\bf N}}

\def \CC {{\bf C}}

\def \im {{\rm {Im\,}}}
\def \G {{\bf G}_m}

\def\Sha{\cyrille X}

\def\ov{\overline}
\def\til{\widetilde}

\def\Si{\Sigma}

\def\Ga{\Gamma}

\def\smallsquare{\vbox{\hrule\hbox{\vrule height 1 ex\kern 1 ex\vrule}\hrule}}
\def\enddem{\hfill \smallsquare\vskip 3mm}
\def \merci {\paragraph{Remerciements. }}
\def \abstract{\paragraph{R\'esum\'e. }}

\DeclareFontFamily{U}{wncy}{}
\DeclareFontShape{U}{wncy}{m}{n}{%
   <5>wncyr5%
   <6>wncyr6%
   <7>wncyr7%
   <8>wncyr8%
   <9>wncyr9%
   <10>wncyr10%
   <11>wncyr10%
   <12>wncyr6%
   <14>wncyr7%
   <17>wncyr8%
   <20>wncyr10%
   <25>wncyr10}{}
\DeclareMathAlphabet{\cyrille}{U}{wncy}{m}{n}
\def\Sha{\cyrille X}

\def\Ga{\Gamma}
\def\Si{\Sigma}
\def \R{{\bf R}}
\def \ss{{\rm ss}}
\def \sc{{\rm sc}}
\def \uu{{\rm u}}
\def \ab{{\rm ab}}
\def \alg{{\rm alg}}
\def \tor{{\rm tor}}
\def \red{{\rm red}}
\def \mult{{\rm mult}}
\def \et {{\rm et}}
\def \fp {{\rm fppf}}
\def \bro {{\rm Br}_0 \,}

\def\Br{{\rm Br\,}}
\def\Pic{{\rm Pic\,}}
\def\Gbar{{\overline{G}}}
\def\Gm{{\mathbf{G}_m}}
\def\SL{{\rm SL}}
\def\labelto#1{\smash{\mathop{\longrightarrow}\limits^{#1}}}

\author{Mikhail Borovoi, Cyril Demarche et David Harari}

\title{Complexes de groupes de type multiplicatif
\\ et groupe de Brauer non ramifi\'e
\\ des espaces homog\`enes}

\begin{document}
\maketitle

\section{Introduction}

Soit $G$ un groupe alg\'ebrique
lin\'eaire connexe lisse d\'efini sur un corps $k$.
Soit $X$ un espace homog\`ene de $G$, tel que le stabilisateur
g\'eom\'etrique $\ov H$ soit {\it de type (ssumult)} ,
c'est-\`a-dire une extension d'un groupe de
type multiplicatif lisse par un groupe connexe, lisse et sans
caract\`eres (par exemple, un groupe lin\'eaire
lisse et connexe sur un corps
parfait est de type (ssumult), ainsi qu'un groupe lin\'eaire
commutatif sur un corps de caract\'eristique z\'ero).
Le but de cet article est d'\'etablir des formules
(obtenues pr\'ec\'edemment dans des cas particuliers
par Colliot-Th\'el\`ene et Kunyavski\u\i\   \cite{ctku1},
\cite{ctku2}) pour le groupe de
Brauer alg\'ebrique $\brr_1 X^c$ (et parfois pour tout le groupe de Brauer
$\br X^c$) d'une compactification
lisse $X^c$ de $X$. Le groupe $\br X^c$
co\"{\i}ncide avec le groupe de Brauer non
ramifi\'e de $X$ si $k$ est de caract\'eristique z\'ero (et le m\^eme
\'enonc\'e vaut en caract\'eristique $p$ exception faite de la torsion
$p$-primaire).

Dans le cas o\`u $k$ est un corps global, on rappelle que ce groupe $\br X^c$
intervient de fa\c con cruciale dans l'\'etude de l'arithm\'etique de l'espace homog\`ene $X$ et de sa compactification $X^c$,
via l'obstruction de Brauer-Manin au principe de Hasse
et \`a l'approximation faible (voir par exemple \cite{skobook}, th\'eor\`eme 5.2.1 (a)).
Par ailleurs, ce groupe est \'egalement utilis\'e pour \'etudier la rationalit\'e de l'espace homog\`ene $X$
sur un corps $k$ quelconque (voir par exemple \cite{chili}). Ces deux exemples illustrent l'int\'er\^et
et la motivation que l'on peut avoir \`a disposer de formules d\'ecrivant le groupe $\br X^c$ sur un corps $k$ quelconque.

Comme dans les travaux ant\'erieurs, les formules obtenues ici pour le groupe $\br X^c$ font intervenir des groupes
"$\Sha_{\omega}$", mais ils sont associ\'es \`a des groupes d'hypercohomologie galoisienne de certains complexes (d\'efinis
et utilis\'es dans \cite{borvh} et \cite{borvh2})
au lieu simplement de groupes de cohomologie galoisienne.

\smallskip

Les g\'en\'eralisations obtenues sont dans plusieurs directions~: tout
d'abord certains r\'esultats valent en caract\'eristique $p$ (notamment
si $k$ est un corps fini ou un corps local),
on ne demande plus que le groupe $G$ soit quasi-trivial mais
seulement de groupe de Picard g\'eom\'etrique trivial, et l'hypoth\`ese
sur le stabilisateur g\'eom\'etrique est plus faible que l'hypoth\`ese
habituelle de connexit\'e. Enfin, sous certaines hypoth\`eses
suppl\'ementaires, on obtient aussi l'\'egalit\'e $\brr_1 X^c=\br X^c$.
Par ailleurs l'utilisation d'un th\'eor\`eme r\'ecent de Gabber (qui
pr\'ecise le th\'eor\`eme de De Jong sur les alt\'erations) permet
de d\'emontrer la plupart des \'enonc\'es sur $\brnr X$ en
caract\'eristique $p$ (en exceptant
la torsion $p$-primaire) sans supposer l'existence d'une compactification
lisse pour $X$.

\smallskip

Une autre nouveaut\'e est que la m\'ethode de d\'emonstration passe
comme dans \cite{ctku1} par le cas des corps finis, mais pour traiter
ceux-ci on est amen\'e \`a passer sur des corps globaux de caract\'eristique
$p$ et \`a utiliser des th\'eor\`emes de dualit\'e arithm\'etique locale
pour obtenir le r\'esultat dans ce contexte. C'est pourquoi on a \'et\'e amen\'e
\`a \'etendre certains r\'esultats classiques de dualit\'e
sur les corps $p$-adiques aux corps locaux de caract\'eristique quelconque,
et aussi
\`a certains complexes de groupes de type multiplicatif au lieu de se limiter
\`a des modules galoisiens (proposition~\ref{perflocal},
proposition~\ref{abelocal}).

\smallskip

Les r\'esultats principaux sont les suivants~:

\smallskip

-Pour un corps $k$ de caract\'eristique z\'ero, on a une formule
qui est valable d\`es que $\ov H$ est de type (ssumult) et $G$ de
groupe de Picard g\'eom\'etrique nul
(th\'eor\`eme~\ref{carzerotheo}). On a de plus des contre-exemples
qui montrent que l'hypoth\`ese que $\ov H$ est de type (ssumult) ne peut
pas \^etre supprim\'ee; en particulier l'hypoth\`ese que le
groupe des composantes connexes de $\ov H$ est commutatif n'est pas
suffisante (proposition~\ref{contrexprop}).

-On montre \'egalement une 
formule analogue sur un corps global de caract\'eristique
$p$ (th\'eor\`eme~\ref{globaltheo}),
et on prouve \'egalement la trivialit\'e du groupe consid\'er\'e sur un corps
fini (th\'eor\`eme~\ref{finitheo}) moyennant quelques hypoth\`eses
de lissit\'e additionnelles.

-Si $H$ est suppos\'e connexe et $k$ de caract\'eristique z\'ero,
on montre que le groupe de Brauer g\'eom\'etrique $\br \ov X^c$ est nul
(th\'eor\`eme~\ref{theo transc})
En particulier la formule du th\'eor\`eme~\ref{carzerotheo} vaut alors en
rempla\c cant $\brr_1 X^c$ par $\br X^c$. On a \'egalement un analogue
partiel en caract\'eristique $p$.

\smallskip

Au passage on donne aussi une d\'emonstration d\'etaill\'ee
d'une formule de compatibilit\'e entre
accouplements (th\'eor\`eme~\ref{compatible}) qui est utile en elle-m\^eme.

\section{Notations et conventions}

Dans ce texte, tous les sch\'emas en groupes sont suppos\'es
localement de type fini. Si $k$ est un corps, on dit simplement
"$k$-groupe" pour d\'esigner un $k$-sch\'ema en groupes lin\'eaire,
i.e. affine et de type fini.
Un
{\it $k$-groupe r\'eductif} (resp.
semi-simple) est un $k$-sch\'ema en groupes affine, lisse, connexe, et
r\'eductif (resp. semi-simple). On note souvent $G^0$ la composante
connexe du neutre d'un $k$-groupe $G$. Si $G$ est un
$k$-groupe connexe lisse sur un corps parfait, on note
$G^{\uu}$ son radical unipotent,
$G^{\red} := G / G^{\uu}$ le quotient r\'eductif, $G^{\ss}$ le sous-groupe d\'eriv\'e de $G^{\red}$ et $G^{\sc}$
le rev\^etement semi-simple simplement connexe de $G^{\sc}$. On dit alors que $G$ est {\it quasi-trivial} si $G^{\tor}:=
G^{\red}/G^{\ss}$ est quasi-trivial (i.e. son groupe des caract\`eres
est un module galoisien de permutation) et $G^{\ss}$ est simplement
connexe.

\smallskip

Pour tout corps $k$, on note $\kbar$ une cl\^oture {\it s\'eparable} de $k$ et
$\Ga_k=\gal(\kbar/k)$.
Si $X$ est un $k$-sch\'ema et $L$ une extension de $k$,
on pose $X_L:=X \times_k L$ et $\ov X=X \times_k \kbar$.
Les groupes de cohomologie \'etale $H^n _{\et}(X,...)$ sont not\'es
simplement $H^n(X,...)$, tandis qu'on note $H^n_{\fp}(X,...)$
s'il s'agit de groupes fppf (les deux co\"{\i}ncident pour des faisceaux
repr\'esent\'es par des sch\'emas en groupes lisses et quasi-projectifs).
Pour simplifier on note aussi souvent $H^n(X,C^{\bullet})$
les groupes d'hypercohomologie $\hyp^n(X,C^{\bullet})$ quand
$C^{\bullet}$ est un complexe
born\'e de faisceaux \'etales sur $X$ (et de m\^eme avec les groupes
fppf). On adopte des notations analogues pour les ensembles de cohomologie
(non ab\'elienne) $H^1(X,...)$. En particulier
si $G$ est un $k$-sch\'ema en groupes
lisse, l'ensemble $H^1(k,G)$ (resp. les groupes $H^i(k,G)$ si $G$ est
commutatif) s'identifie \`a l'ensemble de cohomologie galoisienne
$H^1(\Ga_k,G(\kbar))$  (resp. aux groupes de cohomologie galoisienne
$H^i(\Ga_k,G(\kbar))$).

\smallskip

Par convention, quand on \'ecrit un complexe \`a deux
termes sous la forme $[A \to B]$, cela signifie que $A$ est en degr\'e
$-1$ et $B$ en degr\'e $0$. Si $Z$ est un sch\'ema, on note
${\cal D}(Z)$ (resp. ${\cal D}_{\fp}(Z)$) la cat\'egorie d\'eriv\'ee
born\'ee des faisceaux \'etales (resp. fppf) sur $Z$.
Si $k$ est un corps et $M$ est un $\Ga_k$-module galoisien
(ou encore un complexe born\'e de modules galoisiens), on note
(pour tout $i \geq 1$) $\Sha^i_{\omega, \alg}(k,M)$
(ou $\Sha^i_{\omega ,\alg}(M)$ s'il
n'y a pas d'ambig\"uit\'e sur le corps $k$) le sous-groupe de
$H^i(\Ga_k,M)$ constitu\'e des \'el\'ements dont la restriction \`a
$H^i(C,M)$ est nulle pour tout sous-groupe procyclique $C$ de $\Ga_k$.

\smallskip

Un {\it corps local} est un corps complet pour une valuation discr\`ete
\`a corps r\'esiduel fini~: c'est une extension finie de $\Q_p$ (s'il est
de caract\'eristique $0$) ou d'un corps de s\'eries de Laurent
$\F_q((t))$ sur un corps fini (s'il est de caract\'eristique $>0$). On
utilisera fr\'equemment le fait (d\^u \`a
Kneser \cite{kneser1}, \cite{kneser2} en caract\'eristique
z\'ero et \`a Bruhat-Tits \cite{brtits} en caract\'eristique $>0$) que
pour un groupe semi-simple simplement connexe $G$ sur un corps local $K$,
on a $H^1(K,G)=0$. Par Hilbert 90, on en d\'eduit alors le r\'esultat analogue
si $G$ est quasi-trivial sur un corps local.

\smallskip

Un {\it corps global} est un corps
de nombres (extension finie de $\Q$) ou le corps des fonctions d'une courbe
alg\'ebrique sur un corps fini. Pour toute place $v$ d'un tel corps $K$,
on note $K_v$ le compl\'et\'e de $K$ en $v$.
Si $K$ est un corps global et $M$
un $\Ga_K$-module galoisien (ou un complexe born\'e de modules galoisiens),
on note $\Sha^i_{\omega}(K,M)$ (ou simplement $\Sha^i_{\omega}(M)$)
le sous-groupe de
$H^i(K,M)$ constitu\'e des \'el\'ements dont la restriction \`a
$H^i(K_v,M)$ est nulle pour presque toute place $v$ de $K$.
Si de plus $M$ est un $\Ga_K$-module galoisien de type fini (resp. de
type fini sans torsion), une cons\'equence facile du
th\'eor\`eme de
Cebotarev est que $\Sha^1_{\omega}(M)=\Sha^1_{\omega, \alg}(M)$
(resp. $\Sha^2_{\omega}(M)=\Sha^2_{\omega, \alg}(M)$), \cite{sansuc},
section 2.

\smallskip

Si $A$ est un groupe ab\'elien, et $n \in \NN$, $A[n]$ d\'esigne le sous-groupe de $n$-torsion de $A$.
Si $l$ est un nombre premier, $A\{l\}$ d\'esigne le sous-groupe de torsion $l$-primaire de $A$,
et $A\{l'\}$ le sous-groupe de torsion premi\`ere \`a $l$.

\section{Accouplement entre complexes de groupes de type multiplicatif}

Dans cette section on \'etend l'accouplement $S \times \widehat S \to \G$
entre un groupe de type multiplicatif $S$ et son groupe des caract\`eres
$\widehat S$ \`a des complexes \`a deux termes, et on montre un
r\'esultat de dualit\'e li\'e \`a cet accouplement sur un corps local.

\begin{lem} \label{defpairing}
Soit $Z$ un sch\'ema. Soient $C=[S \to T]$ un complexe de groupes de
type multiplicatif de type fini sur $Z$ et
$\widehat C=[\widehat T \to \widehat S]$
le complexe dual. Alors on a un accouplement canonique
$$ [S \to T] \otimes^{\bf L} [\widehat T \to \widehat S] \to \G [1]$$
dans la cat\'egorie d\'eriv\'ee born\'ee ${\mathcal D}_{\fp}(Z)$
des faisceaux fppf
sur $Z$. De plus, cet accouplement induit
pour tout $n \geq 0$ un isomorphisme canonique
$$H^n _{\fp}(Z,[S \to T]) \to \R^n
\Hom_Z([\widehat T \to \widehat S], \G [1])$$
o\`u $\Hom_Z(...)$ d\'esigne les homomorphismes dans la cat\'egorie
d\'eriv\'ee born\'ee ${\mathcal D}(Z)$ des faisceaux \'etales sur $Z$.
La m\^eme assertion est valable si l'on remplace $\Hom_Z(...)$ par
$\Hom_{{\mathcal D}_{\fp}(Z)}(...)$.

\end{lem}

\dem  On note que par dualit\'e pour les groupes de type multiplicatif,
le complexe $\widehat C$ est le complexe des morphismes
$\fhom^*(C,\G[1])$ correspondant au foncteur "Hom interne" $\fhom(.,\G[1])$
dans la cat\'egorie des complexes de faisceaux fppf sur $Z$.
Alors $\R \fhom(.,\G[1]))$ est le foncteur d\'eriv\'e total de
$\fhom(.,\G[1])$, d'o\`u un morphisme naturel dans ${\cal D}_{\fp}(Z)$~:
$$\widehat C \to \R \fhom(C,\G[1]) \, .$$
Par adjonction entre le produit tensoriel d\'eriv\'e $\otimes^{\bf L}$
et $\R \fhom$, on a un isomorphisme naturel d'adjonction (cf. \cite{weibel},
Th. 10.8.7)~:
$$\Hom_{{\mathcal D}_{\fp}(Z)}(\widehat C \otimes^{\bf L} C, \G[1])=
\Hom_{{\mathcal D}_{\fp}(Z)}(\widehat C, \R \fhom(C, \G[1]))$$
qui donne un morphisme
$$\hat C \otimes^{\bf L} C \to \G[1]$$
fournissant l'accouplement souhait\'e.

\smallskip

Pour la deuxi\`eme assertion, on note pour commencer que le groupe de
droite ne change pas que l'on travaille dans ${\mathcal D}(Z)$ ou dans
${\mathcal D}_{\fp}(Z)$ (via \cite{skobook}, lemme~2.3.7).
L'isomorphisme
se d\'eduit imm\'ediatement par d\'evissage \`a partir des
cas extr\^emes $S=0$ ou $T=0$
(via le lemme des cinq), le r\'esultat \'etant alors connu (loc. cit.).

\enddem

\rem S'il n'y a pas d'hypoth\`ese de lissit\'e sur $S$ et $T$, il
est important de consid\'erer $H^n _{\fp}(X,[S \to T])$, et non le
groupe \'etale correspondant.

\smallskip

Quand $Z=\spec k$ est le spectre d'un corps, le lemme~\ref{defpairing}
donne en particulier un accouplement de cup-produit
\begin{equation} \label{cupdef}
H^1(k,[\widehat T \to \widehat S]) \times H^0_{\fp}(k, [S \to T]) \to
H^2_{\fp}(k,\G)=H^2(k,\G)=\br k
\end{equation}

Si $K$ est un corps local, on a une topologie naturelle sur les
groupes $H^i_{\fp}(K, [S \to T])$ pour $i \geq 0$, qui en fait des
groupes s\'epar\'es (au sens de Hausdorff), localement compacts,
d\'enombrables \`a l'infini, et totalement discontinus. Si $K$ est un
corps local de caract\'eristique z\'ero, on obtient des groupes discrets
pour $i \geq 1$ (et on peut d\'efinir la topologie sur
$H^0_{\fp}(K, [S \to T])$ via le fait que ce groupe est extension d'un
quotient de $H^0(K,T)$ par un sous-groupe du groupe fini $H^1(K,S)$, cf.
\cite{demirn}, section 3). Si $K$ est de caract\'eristique $p >0$,
le groupe $H^1_{\fp}(K,S)$ est en g\'en\'eral infini (il est seulement profini)
si $S$ contient de la
$p$-torsion; la topologie sur les $H^i_{\fp}(K, [S \to T])$ est alors
d\'efinie par le m\^eme proc\'ed\'e que dans \cite{adt}, section III.6.
On a alors l'\'enonc\'e suivant, dans l'esprit du
th\'eor\`eme 3.1 de \cite{demirn}~:

\begin{prop} \label{perflocal}
Soit $K$ un corps local (de caract\'eristique quelconque). Soit
$[S \to T]$ un complexe de $K$-groupes de type multiplicatif. Alors
l'accouplement (\ref{cupdef}) induit un isomorphisme entre le groupe
discret $H^1(K,[\widehat T \to \widehat S])$ et le dual de Pontryagin
$H^0_{\fp}(K,[S \to T])^D$ de $H^0_{\fp}(K,[S \to T])$.
\end{prop}

Ici le dual de Pontryagin $A^D$ d'un groupe topologique ab\'elien $A$ est
le groupe des homomorphismes continus de $A$ dans le groupe discret
$\Q/\Z$.

\dem On proc\`ede par d\'evissage en commen\c cant par le cas o\`u
$T=0$. Comme $S$ est de type multiplicatif, on peut \'ecrire une suite exacte
$$0 \to S \to T_1 \to T_2 \to 0$$
o\`u $T_1$ et $T_2$ sont des tores. On obtient un diagramme commutatif
\`a lignes exactes
{\small
%
$$
\xymatrix@C=0.7cm{
H^1(K,\widehat T_2)\ar[r]\ar[d] &H^1(K,\widehat T_1) \ar[r]\ar[d] &H^1(K,\widehat S)
\ar[r]\ar[d] &H^2(K,\widehat T_2) \ar[r]\ar[d] &H^2(K,\widehat T_1)\ar[d] \\
H^1(K,T_2)^D \ar[r] &H^1(K,T_1)^D \ar[r] &H^1_{\fp}(K,S)^D
\ar[r] &H^0(K,T_2)^D \ar[r] &H^0(K,T_1)^D
}
$$
}
et le r\'esultat r\'esulte alors du
lemme des cinq joint au fait que pour chaque tore $T_i$ les
fl\`eches $H^1(K,\widehat T_i) \to H^1(K,T_i)^D$ et
$H^2(K,\widehat T_i) \to H^0(K,T_i)^D$ induites par (\ref{cupdef}) sont
des isomorphismes (\cite{adt}, Th. III.6.9; noter que pour chaque tore $T_i$,
le groupe $H^1(K,T_i)$ est fini tandis que $H^0(K,T_i)$ et son compl\'et\'e
pour la topologie des sous-groupes ouverts d'indice fini ont m\^eme dual,
\cite{dhsza}, Lemme 2.2).

\smallskip

Dans le cas g\'en\'eral o\`u $S$ et $T$ sont quelconques, on a un
diagramme commutatif
{\small
$$
\xymatrix@C=0.7cm{
H^1(K,\widehat T_2)\ar[r]\ar[d] &H^1(K,\widehat T_1) \ar[r]\ar[d] &H^1(K,\widehat S)
\ar[r]\ar[d] &H^2(K,\widehat T_2) \ar[r]\ar[d] &H^2(K,\widehat T_1)\ar[d] \\
H^1(K,T_2)^D \ar[r] &H^1(K,T_1)^D \ar[r] &H^1_{\fp}(K,S)^D
\ar[r] &H^0(K,T_2)^D \ar[r] &H^0(K,T_1)^D
}
$$
}
et le m\^eme
argument de d\'evissage fonctionne \`a condition de savoir que
pour un groupe de type multiplicatif $S$ quelconque, la fl\`eche
$H^2(K,\widehat S) \to H^0(K,S)^D$ est un isomorphisme. Or ceci est
connu si $S$ n'a pas de $p$-torsion (\cite{adt}, Cor. I.2.3, toujours
avec le fait que $H^0(K,S)$ et son compl\'et\'e ont m\^eme dual),
et le cas g\'en\'eral s'en d\'eduit
imm\'ediatement vu que si $\mu$ est un groupe de type multiplicatif de
torsion $p$-primaire, on a $H^0_{\fp}(K,\mu)=H^2(K,\hat \mu)=0$ et
$H^1(K,\hat \mu) \xrightarrow{\cong} H^1_{\fp}(K, \mu)^D$
(\cite{adt}, Prop. III.6.4 et Prop. III.6.10).

\enddem

\rem On pourrait montrer plus pr\'ecis\'ement qu'on obtient une dualit\'e
entre le groupe discret
$H^1(K,[\widehat T \to \widehat S])$ et le compl\'et\'e (pour la topologie
d\'efinie par les sous-groupes ouverts d'indice fini) de
$H^0_{\fp}(K,[S \to T])$; nous n'aurons pas besoin de ce r\'esultat.

\section{Groupe de Brauer non ramifi\'e, relation avec le
complexe $KD'(X)$} \label{sec:4}

Soit $p : X \to \spec k$ une vari\'et\'e lisse et
g\'eom\'etriquement int\`egre sur un
corps $k$, dont on note $\br X:=H^2(X,\G)$ le groupe de Brauer cohomologique.
Rappelons que $X$ \'etant un sch\'ema int\`egre r\'egulier, la
fl\`eche canonique $\br X \to \br (k(X))$ est injective (o\`u $k(X)$ est le
corps des fonctions de $X$) par le corollaire~II.1.8 de \cite{GrBrauer}.
 Si $X^c$ est une compactification lisse de $X$, alors la fl\`eche compos\'ee
$\br X^c \to \br X\to \br (k(X))$ est injective, donc la
fl\`eche $\br X^c \to \br X$ est injective.

\smallskip

Comme dans \cite{opendesc}, on
consid\`ere le complexe de faisceaux \'etales sur $\spec k$
d\'efini par
$$KD(X)=({\rm UPic} \, X)[1]:=(\tau_{\leq 1} \R p_* {\bf G}_{m,X})[1]$$
et on d\'efinit $KD'(X)$ comme le c\^one du morphisme canonique
$\G[1] \to KD(X)$.

Posons $\bra X=\ker [\br X \to \br \ov X]$. Un argument de suite
spectrale donne des fl\`eches
{\small
$$
H^1(k,\R p_* {\bf G}_{m,X}) \to H^1_{\fp}(k,\R p_* {\bf G}_{m,X}) \to
\br X=H^2(X,\G)=H^2_{\fp}(X,\G)$$
}
qui induisent un isomorphisme $H^1(k,KD(X)) \simeq \bra X$, et une fl\`eche
$$r : \bra X \simeq H^1(k,KD(X)) \to H^1(k,KD'(X))$$ dont le noyau est
$\bro X:=\im [\br k \to \br X]$, fl\`eche qui est de plus surjective si
$H^3(k,\G)=0$ ou si $X(k) \neq \emptyset$ (\cite{borvh}, Prop. 2.18
et \cite{opendesc}, page 7).
Si $m \in X(k)$ est un point rationnel de $X$, on notera
$\brr_{1,m} \,X$ le sous-groupe de $\bra X$ constitu\'e des
\'el\'ements dont l'\'evaluation en $m$ est nulle. Ainsi
$\brr_{1,m} \, X$ s'identifie \`a $\bra X/\br k$.

\begin{prop} \label{facilinclusion}
Soit $X$ une vari\'et\'e lisse, g\'eom\'etriquement int\`egre sur un
corps $k$.
On suppose que $X$ admet une compactification lisse $X^c$. Alors
l'image r\'eciproque du sous-groupe
$\Sha^1_{\omega, \alg}(KD'(X)) \subset H^1(k,KD'(X))$ par $r$
 est contenue dans l'image de $\bra X^c$ dans $\bra X$.
 Le m\^eme
r\'esultat vaut avec $\Sha^1_{\omega}(KD'(X))$ \`a la place de
$\Sha^1_{\omega, \alg}(KD'(X))$ si $k$ est un corps global.
\end{prop}

\dem
D'apr\`es \cite{borvh}, prop. 2.19, on a un diagramme commutatif exacte
\begin{equation}\label{eq:diagBvH}
\xymatrix{
\Br k\ar[r]\ar@{=}[d]  &\bra X^c\ar[r]\ar[d]^{i_{\bra}^*}  &H^1(k,KD'(X^c))\ar[r]\ar@{^{(}->}[d]^{i^*_{KD'}} &H^3(k, \Gm)\ar@{=}[d]\\
\Br k\ar[r]            &\bra X\ar[r]^-r                           &H^1(k,KD'(X))\ar[r]                     &H^3(k, \Gm) \, .
}
\end{equation}
D'apr\`es  \cite{borvh}, cor. 2.16, l'homomorphisme
$i^*_{KD'}$  est injectif et induit un isomorphisme
$$
 \Sha^1_{\omega,\alg}(KD'(X^c))\labelto{\sim} \Sha^1_{\omega,\alg}(KD'(X)).
$$
Soit $\xi\in \bra X$ tel que $r(\xi)\in  \Sha^1_{\omega,\alg}(KD'(X))$, alors
$$r(\xi)\in i^*_{KD'}(\Sha^1_{\omega,\alg}(KD'(X^c)))\subset \im i^*_{KD'},$$
et par chasse au diagramme on obtient que $\xi\in \im i_{\bra}^*$.
L'argument avec $\Sha^1_{\omega}$ \`a la place de $\Sha^1_{\omega,\alg}$ est
identique.
\enddem

\rem Si $k$ est de caract\'eristique z\'ero, l'existence de
$X^c$ est automatique par la r\'esolution des singularit\'es d'Hironaka.
Dans ce cas $\br X^c$ n'est autre que le {\it groupe de Brauer non ramifi\'e}
$\brnr X$ de $X$ (\cite{chili}, section 5). En caract\'eristique $p$,
si $X^c$ existe, on sait juste que pour tout nombre premier $l$
diff\'erent de $p$, on a $(\brnr X)\{ l \}=(\br X^c)\{ l \}$
(\cite{colliotgersten}, proposition~4.2.3).

\smallskip

Montrons maintenant une sorte de r\'eciproque de la proposition \ref{facilinclusion},
sur un corps global, sans supposer l'existence d'une compactification lisse de $X$.
\begin{prop} \label{prop alteration}
Soit $K$ un corps global de caract\'eristique $p \geq 0$, $X$ une $K$-vari\'et\'e lisse et g\'eom\'etriquement int\`egre.
Soit $\alpha \in (\brnr X) \{p'\}$. Alors pour presque toute place $v$ de $K$, $\alpha_v \in \Br(X_v)$
est orthogonal \`a $X(K_v)$ pour l'accouplement
d'\'evaluation $X(K_v) \times \Br(X_v) \to \Br(K_v) \cong \Q / \Z$, o\`u
on a not\'e $X_v:=X \times_K K_v$.
\end{prop}

\dem Supposons tout d'abord
l'existence d'une compactification lisse $X^c$ de $X$ sur $K$ (par
 exemple $K$ de caract\'eristique nulle).
Il existe alors un ensemble fini de places
$\Si$ de $K$ (contenant les \'eventuelles places archim\'ediennes)
tel que
la $K$-vari\'et\'e propre et lisse
$X^c$ s'etende en un sch\'ema propre et
lisse ${\mathcal X}^c$ au-dessus de l'anneau des $\Si$-entiers $\calo_{\Si}$.
Soit $\alpha \in (\brnr X) \{p'\}$, alors $\alpha \in \br X^c$ d'apr\`es la
remarque ci-dessus.
Quitte \`a augmenter $\Si$, on peut donc supposer que
$\alpha$ est dans $\br {\mathcal X}^c$.
Or le groupe de Brauer de l'anneau des entiers $\calo_v$
de tout compl\'et\'e non archim\'edien $K_v$ de $K$ est nul; par
cons\'equent on a $\alpha(P_v)=0$ pour toute $v \not \in \Si$ et tout
point local $P_v \in X(K_v)$ car par propret\'e $P_v$ s'\'etend en une
$\calo_v$-section de ${\mathcal X}^c$.

\smallskip

Traitons maintenant le cas g\'en\'eral o\`u $K$ est de caract\'eristique
$p>0$, ce qui fait que l'existence d'une compactification lisse
n'est pas connue.
On peut supposer que $\alpha \in \brnr(X)\{l\}$, avec $l$ premier distinct de $p$.
Le th\'eor\`eme de compactification de Nagata assure l'existence d'une immersion ouverte $X \to Z$, o\`u $Z$ est une $K$-vari\'et\'e propre. Un
th\'eor\`eme de Gabber (voir \cite{Gab}, expos\'e X, th\'eor\`eme 2.1)
assure qu'il existe une extension finie de corps $K'/K$,
de degr\'e premier \`a $l$, et une $l'$-alt\'eration
$h : Y \to Z_{K'}$, o\`u $Y$ est une $K'$-vari\'et\'e lisse.
En particulier, le morphisme $h$ est propre, surjectif et
g\'en\'eriquement fini de degr\'e premier \`a $l$, et $Y$ est donc propre et
lisse sur $K'$.
On consid\`ere alors le diagramme commutatif suivant, o\`u les deux carr\'es sont des produits fibr\'es :
\begin{equation} \label{diag alteration}
  \xymatrix{
V \ar[r]^i \ar[d]^{h'} & Y \ar[d]^h \\
X_{K'} \ar[r] \ar[d] & Z_{K'} \ar[d] \\
X \ar[r] & Z \, .
}
\end{equation}
On note $\pi$ le morphisme compos\'e $V \to X_{K'} \to X$. Par fonctorialit\'e du groupe de Brauer non ramifi\'e
(voir \cite{chili}, lemme 5.5), on a $\pi^*(\alpha)
\in \brnr(V)\{ l \} = \br(Y) \{l \}$.

Le morphisme $\pi : V \to X$ \'etant g\'en\'eriquement fini de degr\'e $d$ premier \`a $l$, il existe un ouvert affine $U$ de $X$
tel que la restriction $\pi_U : V_U:=V \times_X U \to U$
de ce morphisme soit un morphisme fini et plat sur $U$,
donc dont le degr\'e des fibres est un entier constant $d$ premier \`a $l$.

Soit $S$ un ensemble fini de places de $K$ tel que le diagramme \eqref{diag alteration}
(ou plus exactement son analogue en rempla\c cant $X$, $X_{K'}$, $V$
respectivement par $U$, $U_{K'}$, $V_U$)
s'\'etende en un diagramme de $\spec(\mathcal{O}_{K,S})$-sch\'emas
\begin{displaymath}
\xymatrix{
\mathcal{V_U} \ar[r]^i \ar[d]^{h'} & \mathcal{Y} \ar[d]^h \\
\mathcal{U}_{K'} \ar[r] \ar[d] & \mathcal{Z}_{K'} \ar[d] \\
\mathcal{U} \ar[r] & \mathcal{Z}
}
\end{displaymath}
avec $\mathcal{Y}$ propre lisse sur $\spec(\mathcal{O}_{K',S'})$, o\`u $S'$
d\'esigne l'ensemble des places de $K'$ au-dessus de $S$.
Quitte \`a agrandir $S$, on peut aussi supposer que $\alpha \in \Br(\mathcal U)$
et que $\pi^*(\alpha) \in \br(\mathcal{Y})$ (car $\pi^*(\alpha) \in \br Y$).

\smallskip

Soit $v \notin S$ et $x_v \in U(K_v)$.
Puisque l'on a l'\'egalit\'e $[K':K] = \sum_{w | v} e_w f_w$ (voir par exemple \cite{rosen}, th\'eor\`eme 7.6),
le fait que $l$ ne divise pas $[K':K]$ assure qu'il existe une place $w$ de $K'$ au-dessus de $v$
telle que $[K'_w:K_v] = e_w f_w$ soit premier \`a $l$.

Puisque $\pi_U : V_U \to U_{K'}$
est fini de degr\'e premier \`a $l$,
il existe une extension finie de corps $L/K'_w$ et un point $y \in X(L)$
tels que $x_v = \pi_U(y)$, avec $[L:K'_w]$ (et donc aussi $[L:K_v]$)
premier \`a $l$.

On a alors $\textup{Res}_{L/K_v}(\alpha(x_v)) = ((\pi_U)^*(\alpha)) (y)$ dans
$\Br(L)$. Or on dispose du diagramme commutatif suivant :
\begin{displaymath}
\xymatrix{
\spec(L) \ar[r]^y \ar[d] & \mathcal{Y} \ar[d] \\
\spec(\mathcal{O}_L) \ar[r] & \spec(\mathcal{O}_{K',S'}) \, ,
}
\end{displaymath}
donc le crit\`ere valuatif de propret\'e assure  
l'existence d'une factorisation $$y : \spec(\mathcal{O}_L) \to \mathcal{Y} .$$
Comme $(\pi_U)^*(\alpha) \in \br \mathcal{Y}$, on en d\'eduit que
$((\pi_U)^*(\alpha)) (y) \in \br(\mathcal{O}_L) = 0$.
Donc finalement, on a montr\'e que $\textup{Res}_{L/K_v}(\alpha(x_v)) = 0$
dans $\br(L)$. Comme $[L:K_v]$ est premier \`a $l$ et que
la restriction $\br K_v \to \br L$ correspond \`a la multiplication par
$[L:K_v]$ (qui est premier \`a $l$) dans $\Q/\Z$, on obtient que
$\alpha(x_v)=0$ dans $\br(K_v)$. Par cons\'equent, l'\'el\'ement $\alpha_v \in \br(X_v)$
est orthogonal \`a $U(K_v)$ pour toute place $v \notin S$. Enfin, puisque $U$ est un ouvert non vide
de la vari\'et\'e int\`egre $X$, le th\'eor\`eme
des fonctions implicites donne que
$U(K_v)$ est dense dans $X(K_v)$ pour la topologie $v$-adique. Or l'accouplement $\br(X_v) \times X(K_v) \to \Br(K_v)$
est continu pour cette topologie, donc $\alpha_v$ est orthogonal \`a $X(K_v)$ tout entier, pour tout $v \notin S$.
Cela conclut la preuve.

\enddem

\section{Le groupe de Brauer g\'eom\'etrique}\label{sec:transc}

Dans cette section, on montre que dans le cas des espaces homog\`enes
\`a stabilisateurs connexes sur un corps alg\'ebriquement clos, le groupe
de Brauer d'un mod\`ele projectif lisse est trivial.
Plus pr\'ecis\'ement le th\'eor\`eme principal de cette section est le suivant :
\begin{theo} \label{theo transc}
Soit $k$ un corps s\'eparablement clos. 

Soit $G$ un 
$k$-groupe lin\'eaire connexe lisse et soit 
$H$ un $k$-sous-groupe connexe lisse de $G$.
On note $X := G/H$ et on d\'esigne par 
$X^c$ une compactification lisse de $X$. 
On suppose l'une des conditions suivantes
\begin{itemize}
    \item $k$ est de caract\'eristique nulle.
    \item $k$ est alg\'ebriquement clos (i.e. $k$ est parfait) et $H$ est produit semi-direct de $H^{\uu}$ par $H^{\red}$.
    \item $H$ et $G$ sont r\'eductifs.
\end{itemize}
Alors $\Br(X^c)\{l\} = 0$, pour tout $l$ premier \`a la caract\'eristique de $k$.
\end{theo}

\begin{rem}
Ce r\'esultat \'etend le r\'esultat classique de Bogomolov (voir \cite{bog}, th\'eor\`eme 2.4),
qui concerne le cas o\`u $k = \CC$ et $G$ est semi-simple simplement connexe.
Cela r\'epond en particulier par l'affirmative \`a la remarque 9.14 de \cite{chili}, ainsi qu'\`a la question qui suit le th\'eor\`eme 1.4
de \cite{ctku2}, o\`u les auteurs demandent si le r\'esultat de Bogomolov s'\'etend (en caract\'eristique nulle)
au cas o\`u $G$ est un groupe lin\'eaire connexe quelconque.
\end{rem}

\smallskip

On commence par un lemme g\'en\'eral~:

\begin{lem} \label{alglem}
Soit $X$ une vari\'et\'e lisse et g\'eom\'etriquement int\`egre sur un corps
$k$. Soient $l$ un nombre premier diff\'erent de ${\rm Car } \, k$ et
$\alpha \in (\br X)\{ l \}$. On suppose que $\alpha \in \brnr X$ (ou encore
$\alpha \in \br X^c$ si on dispose d'une compactification lisse $X^c$ de $X$)
et qu'il existe une extension de corps $K$ de $k$ telle que la restriction
$\alpha_K \in \br X_K$ soit nulle. Alors $\alpha$ est alg\'ebrique (i.e.
la restriction de $\alpha$ \`a $\br \ov X$ est nulle).
\end{lem}

\dem Choisissons des cl\^otures s\'eparables compatibles $\kbar$ de $k$ et
$\ov K$ de $K$. L'hypoth\`ese implique que la restriction de $\alpha$
\`a $\brnr (X_{\ov K})\{ l \}$ est nulle.

Supposons d'abord que $X$ admet une compactification lisse $X^c$. Alors la restriction
$\brnr (X_{\kbar})\{ l \} \to \brnr (X_{\ov K})\{ l \}$ est un isomorphisme
via \cite{colliotgersten}, proposition~4.2.3 et la rigidit\'e
de la cohomologie non ramifi\'ee \`a coefficients finis (loc. cit.,
th\'eor\`eme~4.1.1), d'o\`u le lemme.

Dans le cas g\'en\'eral, on dispose du th\'eor\`eme de Gabber : comme dans la preuve de la proposition \ref{prop alteration},
il existe une extension finie $k'/k$ de degr\'e premier \`a $l$, que l'on peut supposer contenue dans $K$,
et un diagramme cart\'esien analogue \`a \eqref{diag alteration} :
\begin{equation}
  \xymatrix{
V \ar[r]^i \ar[d]^{h'} & Y \ar[d]^h \\
X_{k'} \ar[r] \ar[d] & Z_{k'} \ar[d] \\
X \ar[r] & Z \, ,
}
\end{equation}
avec $Y$ propre et lisse sur $\spec k'$, $V \to Y$ immersion ouverte et $V \to X_{k'}$ g\'en\'eriquement fini de degr\'e premier \`a $l$.

L'argument de rigidit\'e pr\'ec\'edent, appliqu\'e \`a $V$ (qui admet une compactification lisse)
assure que la restriction de $\alpha$ \`a $\br(V_{\overline{k'}})$ est nulle.
Or le morphisme $V_{\ov{k'}} \to X_{\ov{k'}}$ est g\'en\'eriquement fini, de degr\'e premier \`a $l$,
donc un argument de restriction-corestriction au niveau de l'extension de corps de fonctions
assure que la restriction de $\alpha$ \`a $\br(X_{\ov{k'}})$ est nulle
(on utilise aussi ici que si $F/E$ est une extension de corps de caract\'eristique $p$ purement ins\'eparable,
et $W$ une $E$-vari\'et\'e, alors la restriction $\br(W)\{p'\} \to \br(W_F)\{p'\}$ est un isomorphisme : voir \cite{milne}, remarque II.3.17).

Par cons\'equent, il existe une extension finie $L/k$ telle que la restriction de $\alpha$ \`a $\br(X_L)$ soit nulle.
On d\'ecompose alors l'extension $L/k$ en une extension s\'eparable $E/K$ et une extension purement ins\'eparable $L/E$.
La remarque II.3.17 de \cite{milne} assure alors que la restriction $\br(X_E)\{p'\} \to \br(X_L)\{p'\}$ est un isomorphisme,
donc la restriction de $\alpha$ \`a $\br(X_E)$ est nulle, donc $\alpha$ est alg\'ebrique.
\enddem

\paragraph{Preuve du th\'eor\`eme~\ref{theo transc}.}
La preuve de ce r\'esultat se fait en plusieurs \'etapes.

\paragraph{\'Etape 1 :}cas de la cl\^oture s\'eparable d'un corps global.


\begin{lem} \label{lem transc global}
Soit $k$ un corps global. Soit $l$ premier \`a $\textup{Car}(k)$.
\begin{itemize}
  \item Soient $H$ et $G$ deux $k$-groupes r\'eductifs, avec $H \subset G$. On pose $X := G/H$. Soit $X^c$ une compactification lisse de $X$.

Alors $\Br(X^c)\{l\} = \Br_1(X^c)\{l\}$.
\item Si $H$, $G$, $X$ et $X^c$ sont d\'efinis sur $\overline{k}$, alors $\Br(X^c)\{l\}=0$.
\end{itemize}
\end{lem}

\begin{rem}
Gr\^ace \`a la proposition \ref{prop alteration} et au lemme \ref{alglem},
ce lemme reste valable sans supposer l'existence de la compactification lisse $X^c$ de $X$,
en rempla\c cant $\Br(X^c)$ par $\brnr(X)$.
Noter aussi qu'en particulier $H$ et $G$ sont suppos\'es lisses.
\end{rem}

\dem
Tout d'abord, il est clair que le second point du lemme est une cons\'equence
du premier. En effet, si tout est d\'efini sur $\ov k$, et si
$\alpha \in \br(X^c)\{l\}$, alors il existe une extension finie
s\'eparable $K/k$ telle que $H$, $G$, $X$, $X^c$ et $\alpha$ soient
d\'efinis sur $K$. Or $K$ est un corps global, donc $\alpha$ vu dans
le groupe de Brauer de la $K$-vari\'et\'e $X^c$ est alg\'ebrique
par le premier point, donc $\alpha$ est nul dans le groupe de Brauer de
la $\ov k$-vari\'et\'e $X^c$.

Montrons maintenant le premier point. D'apr\`es \cite{ctku2},
lemme~1.5, on peut supposer que le groupe
$G$ est quasi-trivial. Soit $\alpha \in \Br_e(X^c)\{l\}$.
On souhaite montrer que $\alpha$ est alg\'ebrique, c'est-\`a-dire
que sa restriction \`a $\br \ov X^c$ est nulle.

On note $\pi : G \to X$ le morphisme quotient.

Quitte \`a remplacer $k$ par une extension finie s\'eparable,
on peut supposer que les groupes r\'eductifs $H$ et $G$ sont $k$-d\'eploy\'es.

Alors $\pi^*(\alpha) = 0$ dans $\brnr_e(G)$ (puisque $\brnr(G)\{l\} = \br(k)\{l\}$, vu que $G$ est r\'eductif d\'eploy\'e, donc rationnel :
voir \cite{borel}, corollaires 14.14 et 18.8).

Par cons\'equent, la suite exacte du th\'eor\`eme 2.4 de \cite{bordem} (valable en caract\'eristique positive,
avec la m\^eme preuve, en supposant $H$ et $G$ r\'eductifs) assure qu'il existe $p \in \Pic(H)$ tel que $\delta(p) = \alpha$,
o\`u $\delta : \Pic(H) \to \Br(X)$ est le morphisme d\'efini par Colliot-Th\'el\`ene et Xu (\cite{ctx}, section 2, p.314).

Puisque $\alpha \in \Br(X)$ est
non ramifi\'e (i.e. appartient \`a $\br X^c$)
et que le groupe de Brauer de l'anneau des entiers $\calo_v$
est nul, on obtient que pour presque toute place $v$ de $k$, la localisation
$\alpha_v \in \Br(X_v):= \br (X \times_k k_v)$
est orthogonale \`a $X(k_v)$, pour l'accouplement canonique d'\'evaluation (voir proposition \ref{prop alteration}).

On consid\`ere le diagramme commutatif suivant :
\begin{equation}
\label{diag BM}
\xymatrix{
X(k_v) \ar[d] & \times & \Br(X_v) \ar[r] & \Br(k_v) \ar[d]^= \\
H^1(k_v,H) & \times & \Pic(H_v) \ar[r] \ar[u]^{\delta} & \Br(k_v) \, .
}
\end{equation}
Dans ce diagramme, le premier accouplement horizontal est l'accouplement
d'\'evaluation usuel.
Pour la d\'efinition de la seconde ligne et la commutativit\'e du diagramme, voir \cite{ctx}, proposition 2.9.

Or $G$ est quasi-trivial, donc $H^1(k_v,G)=1$ pour toute place finie
$v$ de $k$. Donc l'application caract\'eristique $X(k_v) \to H^1(k_v, H)$ est surjective.
Donc le diagramme \eqref{diag BM} assure que, pour presque toute place $v$, $p_v \in \pic
H_v$ est orthogonal \`a $H^1(k_v,H)$.

\smallskip

Soit $T_H$ un tore maximal de $H$
et $T_{H^{\sc}}$ son image r\'eciproque par le morphisme canonique
$H^{\sc} \to H$. Notons
$C_H := [T_{H^{\sc}} \to T_H]$ et
$\widehat{C}_H := [\widehat{T}_H \to \widehat{T}_{H^{\sc}}]$.
Rappelons qu'on dispose alors
du groupe $H^1_{\textup{ab}}(k_v,H):=\mathbf{H}^1(k_v, C_H)$ et d'un
morphisme d'ab\'elianisation $\ab^1_H : H^1(k_v,H) \to H^1_{\textup{ab}}(k_v, H)$
qui est surjectif pour toute place $v$ de $k$ (\cite{borams},
th\'eor\`eme~5.4 en caract\'eristique z\'ero. En caract\'eristique positive,
le th\'eor\`eme 5.5 et l'exemple 5.4(i) de \cite{gonz} assurent le r\'esultat).

On a un diagramme commutatif :
\begin{displaymath}
\xymatrix{
H^1(k_v,H) \ar[d]^{\textup{ab}^1_H} & \times & \Pic\!(H_v) \ar[r] & \Brr(k_v) \ar@{=}[d] \\
H^1_{\textup{ab}}(k_v,H) & \times & H^1(k_v, \widehat{C}_H) \ar[u]^{\cong} \ar[r]^-{\cup} & \Brr(k_v) \, .
}
\end{displaymath}

Or le morphisme $\mathbf{H}^1(k_v, \widehat{C}_H) \to \Pic\!(H_v)$
est un isomorphisme (voir \cite{borvh}, corollaire 5) et
l'application $\textup{ab}^1_H$ est surjective,
donc l'image $q_v$ de $p_v$ dans $\mathbf{H}^1(k_v, \widehat{C}_H)$
est orthogonale \`a $H^1_{\textup{ab}}(k_v, H)$ pour le cup-produit.
Or l'accouplement $H^1_{\textup{ab}}(k_v,H) \times
\mathbf{H}^1(k_v, \widehat{C}_H) \to \Brr(k_v)$
est une dualit\'e parfaite de groupes
finis (voir \cite{demirn}, th\'eor\`eme 3.1 pour la caract\'eristique
nulle; la preuve en caract\'eristique positive est similaire),
donc $q_v = 0$, donc $p_v = 0$.

Finalement, on a montr\'e que $p_v = 0$ pour presque toute place $v$, donc
on a aussi $\alpha_v=0$ pour presque toute place $v$. Ceci implique en
particulier que $\alpha$ est alg\'ebrique via le lemme~\ref{alglem}.

\enddem

\rem Sur un corps $p$-adique, la dualit\'e entre
$H^1(k_v,H)$ et $\pic H_v$ (variante de r\'esultats de Kottwitz)
avait d\'ej\`a \'et\'e obtenue
par Colliot-Th\'el\`ene (\cite{colflasque}, th\'eor\`eme~9.1) en utilisant
les r\'esolutions flasques.

\paragraph{\'Etape 2 :}cas d'un corps alg\'ebriquement clos de caract\'eristique nulle.

Dans cette partie, on d\'emontre le th\'eor\`eme \ref{theo transc}
dans le cas o\`u $k$ est un corps alg\'ebriquement clos de caract\'eristique nulle.

On suppose d'abord que $H$ et $G$ sont r\'eductifs. On voit $\overline{\Q}$ comme un sous-corps de $k$.
La th\'eorie des groupes r\'eductifs assure qu'il existe
des $\overline{\Q}$-groupes r\'eductifs $H'$ et $G'$ tels que $H'_k = H$ et $G'_k = G$.
Par un th\'eor\`eme de Vingberg r\'edig\'e par Margaux
(voir th\'eor\`eme 1.1 de \cite{margaux}), appliqu\'e \`a l'inclusion (d\'efinie sur $k$) $i : H \to G$,
il existe une inclusion $i' : H' \to G'$ d\'efinie sur $\overline{\Q}$ et $g \in G(k)$ tels que $i = \textup{int}(g) \circ i'_k$.
 On dispose alors d'une $\overline{\Q}$-vari\'et\'e $X' := G' / H'$.

Le $k$-morphisme $\textup{int}(g) : G \to G$ induit alors un $k$-isomorphisme $\varphi : X'_k \to X$.
On en d\'eduit donc un isomorphisme de groupes $\Br(X^c) \cong \Br({X'}^c_k)$.

Or la rigidit\'e de la cohomologie non ramifi\'ee \`a coefficients
finis assure que $\Br({X'}^c) = \Br({X'}^c_k)$
(\cite{colliotgersten}, th\'eor\`eme 4.4.1).

Le lemme \ref{lem transc global} de l'\'etape 1 assure que $\Br({X'}^c) = 0$, d'o\`u finalement $\Br(X^c) = 0$.

Supposons maintenant seulement $H$ r\'eductif. On d\'efinit
$Y := X/G^{\uu}=G/(H . G^{\uu})$.
Alors $Y$ est un espace homog\`ene de $G^{\red}$ \`a stabilisateur r\'eductif, donc par la preuve pr\'ec\'edente,
on a $\Br(Y^c) = 0$ si $Y^c$ est une compactification
lisse de $Y$. D'autre part
la fibre g\'en\'erique du morphisme $X \to Y$ est un espace homog\`ene $W$
de $G^u$, elle admet donc un point rationnel (\cite{bocrelle}, Lemme~3.2)
sur le corps des fonctions $K$ de $Y$. En particulier on a un $K$-morphisme
$G^u_K \to W$ qui est un torseur sous un $K$-groupe unipotent, donc un
torseur trivial.
Ceci implique que $W$ est $K$-rationnel (puisque
c'est le cas de tout groupe unipotent en caract\'eristique z\'ero).
Finalement la vari\'et\'e $X$ est stablement $k$-birationnelle \`a $Y$,
donc $\Br(X^c) \cong \Br(Y^c)$, d'o\`u $\Br(X^c)=0$.

\smallskip

On ne suppose plus maintenant ni $H$, ni $G$ r\'eductifs.
Le th\'eor\`eme de Mostow assure l'existence d'une d\'ecomposition en
produit semi-direct $H = H^{\uu} \rtimes H^{\red}$.
On consid\`ere alors le quotient $Z := G / H^{\red}$.
On dispose d'un morphisme naturel $Z \to X$, et comme ci-dessus
$Z$ est stablement $k$-birationnelle \`a $X$, d'o\`u
$\Br(X^c) \cong \Br(Z^c)$. Or $\Br(Z^c) = 0$ par le cas pr\'ec\'edent
($Z$ est un espace homog\`ene de $G$ \`a stabilisateur r\'eductif
$H^{\red}$), donc $\Br(X^c) = 0$, ce qui conclut la preuve du
th\'eor\`eme \ref{theo transc} dans le cas de caract\'eristique nulle.

\paragraph{\'Etape 3 :}cas d'un corps s\'eparablement clos de caract\'eristique positive.

En caract\'eristique positive, le th\'eor\`eme de Vingberg-Margaux n'est pas v\'erifi\'e pour des groupes r\'eductifs semi-simples :
il vaut pour les groupes lin\'eairement r\'eductifs (voir \cite{margaux}, remarque 5.2 pour un contre-exemple).
Par cons\'equent, on ne peut adapter la m\'ethode de l'\'etape 2 pour d\'eduire le th\'eor\`eme du cas des corps globaux.
On va donc utiliser une autre m\'ethode,
consistant \`a se ramener aux corps finis.

Notons d'abord qu'on peut supposer $G$ et $H$
r\'eductifs, car sinon les hypoth\`eses suppl\'ementaires ($k$ parfait
et $H$ produit semi-direct de $H^u$ par $H^{\red}$) permettent de
se ramener \`a cette situation comme \`a la fin de l'\'etape 2.

On commence par les corps finis :
\begin{lem} \label{lem fini}
Si $k$ est la cl\^oture alg\'ebrique d'un corps fini $\F$, alors on a
$\brnr(X)\{l\}=0$, pour tout $l \neq \textup{Car}(\F)$.
\end{lem}

\dem
On consid\`ere l'extension de corps $k \subset K:=\overline{\F(t)}$. Par le lemme \ref{lem transc global}
et la remarque qui le suit, on sait que $\brnr(X_K)\{l\}=0$.
Alors le lemme \ref{alglem} assure que tous les \'el\'ements de $\brnr(X)\{l\}$ sont alg\'ebriques.
Or $k$ est s\'eparablement clos, donc $\brnr(X)\{l\} = 0$.
\enddem

D\'emontrons maintenant le th\'eor\`eme \ref{theo transc} dans le cas g\'en\'eral,
o\`u $k$ est un corps s\'eparablement clos de caract\'eristique $p>0$.
Soit $\alpha \in \Br(X^c)\{l\}$.
Alors il existe une sous-extension de type fini $L/\F_p$ de $k/\F_p$, telle que $H$, $G$, $X$ et $\alpha$ soient d\'efinis sur $L$.

En utilisant la suite exacte
$$0 \to  \Pic(H_K) \xrightarrow{\delta} \Br(X_K) \xrightarrow{\pi^*} \Br(G_K)$$
pour toute extension $K/L$, ainsi que le fait que
$\Br_{\textup{nr}}(G_{\overline{L}})\{l\}=0$ (puisque
$G_{\overline{L}}$ est rationnel car d\'eploy\'e), on sait qu'il existe une
extension finie s\'eparable $K/L$ et $p \in \Pic(H_K)$ tels que
l'image $\alpha_K$ de $\alpha$ dans $\br X_K$ v\'erifie
$\alpha_K = \delta(p)$, avec de plus la propri\'et\'e que
$H^{\tor}$ soit $K$-d\'eploy\'e. On peut en outre supposer que le groupe fondamental $\mu_H$ de $H^{\ss}$ est d\'eploy\'e par $K$.

Pour montrer le th\'eor\`eme, il suffit
de montrer que $\alpha_K$=0.

\smallskip

Puisque $\Pic(H^{\tor}_K) = 0$, le morphisme $\Pic(H_K) \to \Pic(H^{\ss}_K)$ est injectif.

On sait que $p$ est repr\'esentable par une extension centrale de $K$-groupes alg\'ebriques
\begin{equation}
\label{sec}
1 \to {\G}_K \to H_0 \to H_K \to 1 \, .
\end{equation}
On note $H_1$ le pull-back de cette extension par $H^{\ss}_K \to H_K$.

La trivialit\'e de $\Pic(H^{\sc}_K)$ assure que le pull-back de cette suite exacte par $H^{\sc}_K \to H^{\ss}_K$ est une suite scind\'ee.

Il existe alors une $\F_p$-alg\`ebre int\`egre de type fini $B$ de corps des fractions $K$, des $B$-sch\'emas en groupes
$\mathcal{H}$, $\mathcal{G}$, $\mathcal{H}^{\ss}$, deux $B$-sch\'emas $\mathcal{X}$ et $\mathcal{X}^c$ tels que
\begin{itemize}
    \item les fibres g\'en\'eriques de $\mathcal{H}$, $\mathcal{G}$, $\mathcal{H}^{\ss}$, $\mathcal{X}$ et $\mathcal{X}^c$
sont respectivement $H_K$, $G_K$, $H^{\ss}_K$, $X_K$ et $X^c_K$.
    \item $\mathcal{G}$ est un $B$-sch\'ema en groupes quasi-trivial
(extension d'un tore quasi-trivial par un sch\'ema en groupes semi-simple
simplement connexe; en particulier le groupe de Picard des fibres g\'eom\'etriques
de $\mathcal{G}$ est nul).
    \item $\mathcal{H}$ est un $B$-sch\'ema en groupes r\'eductif, et on a une suite exacte centrale de $B$-sch\'emas en groupes
$$1 \to \mu_{H,B} \to \mathcal{H}^{\sc} \to \mathcal{H}^{\ss} \to 1$$
avec $\mathcal{H}^{\sc}$ semi-simple simplement connexe et $\mu_{H,B}$ fini de type multiplicatif et d\'eploy\'e.
    \item la suite \eqref{sec} s'\'etend en une suite exacte centrale sur $B$
$$1 \to {\G}_B \to \mathcal{H}_0 \to \mathcal{H} \to 1 \, ,$$
dont la classe $p_B$ dans $\pic \mathcal{H}$ a pour fibre g\'en\'erique $p \in
\pic H_K$. On peut de m\^eme \'etendre l'extension associ\'ee \`a $H_1$, ce qui
donne une classe $p'_B \in \pic \mathcal{H^{\ss}}$ (image de $p_B$), dont
la fibre g\'en\'erique $p'$ est l'image de $p$ dans $\pic H^{\ss}_K$. On peut
\'egalement supposer que la suite d\'eduite par pull-back via
$\mathcal{H}^{\sc} \to \mathcal{H}^{\ss}$ est scind\'ee, ce qui implique que
l'image de $p'_B$ dans $\pic \mathcal{H}^{\sc}$ est nulle.

    \item $H_K \to G_K$ s'\'etend en un sous-sch\'ema en groupes ferm\'e
 $\mathcal{H} \to \mathcal{G}$ et $\mathcal{X}$ s'identifie au quotient $\mathcal{G} / \mathcal{H}$.
    \item $x \in X(k)$ s'\'etend en $x \in \mathcal{X}(B)$.
    \item $\mathcal{X}^c$ est un $B$-sch\'ema projectif et lisse.
    \item $X \to X^c$ se prolonge en $\mathcal{X} \to \mathcal{X}^c$.
    \item pour tout $s \in \spec(B)$, la fibre $\mathcal{X}^c_s$ est une compactification lisse de $\mathcal{X}_s$ sur $k(s)$.
    \item la classe $\alpha_K \in \Br(X^c_K)$ s'\'etend en une classe $A \in \Br(\mathcal{X}^c)$,
qui co\"{\i}ncide avec l'image de l'extension $\mathcal{H}_1$ (vue dans $\Pic(\mathcal{H})$) via le morphisme
 $\Pic(\mathcal{H}) \to \Br(\mathcal{X})$ (construit comme dans le cas des corps, via le
$\mathcal{H}$-torseur $\mathcal{G} \to \mathcal{X}$ et la suite exacte $1 \to {\G}_B \to \mathcal{H}_0 \to \mathcal{H} \to 1$).
\end{itemize}

Fixons maintenant un point ferm\'e $s \in \spec(B)$. Par construction, le corps r\'esiduel de $s$, not\'e $k(s)$, est un corps fini.

Les propri\'et\'es des mod\`eles qu'on a consid\'er\'es assurent
que la $k(s)$-compactification lisse $\mathcal{X}_s \to \mathcal{X}^c_s$ v\'erifie les hypoth\`eses du lemme \ref{lem fini}.
La classe $A_s \in \Br(\mathcal{X}^c_s)$ est l'image de la classe $p_s$ de ${\mathcal{H}_0}_s$
(vue dans $\Pic(\mathcal{H}_s)$) par le morphisme $\delta_s : \Pic(\mathcal{H}_s) \to \Br(\mathcal{X}_s)$.
 Notons $p'_s$ l'image de $p_s$ dans $\Pic(\mathcal{H}^{\ss}_s)$, alors $p'_s$ est
aussi la fibre en $s$ de $p'_B \in \pic (\mathcal{H}^{\ss})$.
Comme $\pic \ov H_s $ s'injecte dans $\br  \ov X_s $ (car $\pic \ov G_s=0$)
et $A_s$ est alg\'ebrique (lemme \ref{lem fini}), on obtient que
l'image de $p_s$ dans $\pic \ov H_s$ est nulle. Mais $\pic H_s^{\ss}$
s'injecte dans $\pic \ov H_s^{\ss}$ puisque $H_s^{\ss}$ est semi-simple;
on en d\'eduit que $p'_s = 0$.

\smallskip

Or on dispose d'une suite exacte de groupes
$$0 \to H^0(B, \widehat{\mu_{H,B}}) \to \textup{Ext}^c_B(\mathcal{H}^{\ss}, {\G}_B) \to \textup{Ext}^c_B(\mathcal{H}^{\sc}, {\G}_B) \, ,$$
o\`u $\ext^c_B(...)$ d\'esigne le groupe ab\'elien des classes d'extensions centrales de $B$-faisceaux en groupes (voir \cite{ctx}, p.313).

Comme l'image de $p'_B$ dans $\pic {\mathcal{H}^{\sc}}$ est nulle,
la classe $p'_B$ provient d'un \'el\'ement
$\beta \in H^0(B, \widehat{\mu_{H,B}})$, dont l'image $\beta_s$ dans
$H^0(k(s), \widehat{\mu_{H,s}})$ est nulle (car $p'_s=0$). Or $\widehat{\mu_{H,B}}$ est un $B$-groupe constant car
$\mu_{H,B}$ est fini de type multiplicatif d\'eploy\'e; ceci implique que
le morphisme
$$H^0(B, \widehat{\mu_{H,B}}) \to H^0(k(s), \widehat{\mu_{H,s}})$$
est injectif.
Alors $\beta_s = 0$ implique que $\beta = 0$ dans
$H^0(B, \widehat{\mu_{H,B}})$ ce qui implique
$p'_B = 0$. Ainsi la fibre g\'en\'erique $p'$ de $p'_B$ est nulle dans
$\Pic(H^{\ss}_K)$. Enfin, le morphisme
$\Pic(H_K) \to \Pic(H^{\ss}_K)$ est injectif, d'o\`u $p=0$ dans
$\Pic(H_K)$, donc $\alpha_K = 0$ dans $\Br(X_K)$, ce qui ach\`eve
la preuve du th\'eor\`eme \ref{theo transc}.
\enddem

\section{Une formule de compatibilit\'e}\label{sec:compat}

On consid\`ere un $k$-groupe lisse et connexe $G$ (suppos\'e r\'eductif
pour $k$ non parfait) et un $k$-sous-groupe alg\'ebrique
$H$ de $G$. Rappelons qu'alors le $k$-groupe semi-simple
$G^{\ss}$ est simplement connexe si et seulement si $\pic \ov G=0$
via \cite{sansuc}, lemme 6.9(iii) et remarque 6.11.3 joints au fait
que le noyau du
rev\^etement universel $G^{\rm sc} \to G^{\ss}$ est un $k$-groupe fini de
type multiplicatif \cite{harder}.

\begin{defi}
{\rm Soit $H$ un $k$-groupe alg\'ebrique. On dit que $H$ est
{\it sans caract\`eres} si le seul morphisme de $\kbar$-groupes de $\ov H$ dans $\G$ est le morphisme constant.
On dit que $H$ est
{\it de type (ssumult)} s'il existe une suite exacte de $k$-groupes
$$1 \to L \to H \to S \to 1$$
v\'erifiant~: $S$ est lisse de type multiplicatif et $L$ est lisse,
connexe, sans caract\`eres.
}
\end{defi}

Noter en particulier qu'un groupe de type (ssumult) est lisse.
Si $k$ est de caract\'eristique z\'ero,
tout $k$-groupe connexe ou encore tout $k$-groupe commutatif
est de type (ssumult). Si $k$ est parfait de caract\'eristique $p >0$,
tout $k$-groupe lisse et connexe, ou encore tout $k$-groupe lisse commutatif
avec $H/H^0$ sans $p$-torsion, est de type (ssumult); en effet ceci
r\'esulte
de la structure des groupes alg\'ebriques sur un corps parfait. Plus
g\'en\'eralement, si $k$ est parfait, un $k$-groupe lisse $H$ est
de type (ssumult) si et seulement s'il v\'erifie la condition de
\cite{BCS}, section 3
que le noyau de $H \to H^{\mult}$ est connexe sans caract\`eres,
o\`u $H^{\mult}$
est le quotient maximal de $H$ qui est de type multiplicatif.
La situation est plus compliqu\'ee
sur un corps imparfait, parce qu'on n'a pas en g\'en\'eral de radical
unipotent d\'efini sur $k$.

\smallskip

Soit $X:=G/H$ un espace homog\`ene de $G$ poss\'edant un $k$-point.
D'apr\`es \cite{giraud}, Prop. III.3.1.1, on a pour tout $k$-sch\'ema
$S$ une bijection fonctorielle $\tau_S : X(S) \to H^0(S,[H \to G])$
qui induit par changement de base, pour tout $k$-point $s$ de $S$,
une bijection $X(k) \to H^0(k,[H \to G])$. Pour $S=X$, la classe
$\tau_X({\rm id})$ correspondant \`a l'identit\'e de $X$ donne
une classe $[[X]] \in H^0_{\fp}(X,[H \to G])$
induite par la trivialisation du
$X$-torseur (sous $G$) $G \wedge^H G$ donn\'ee par le neutre $e \in G(k)$.
On obtient ainsi une bijection
$$\tau : X(k) \to H^0_{\fp}(k,[H \to G]) \quad x \mapsto [[X]](x)$$
induite par l'\'evaluation.

\begin{theo} \label{compatible}
Soit $G$ un $k$-groupe lisse et connexe, r\'eductif pour $k$ non parfait,
v\'erifiant $\pic \ov G=0$.
On pose $T=G^{\tor}$.
Soit $H$ un $k$-sous-groupe lisse de $G$, suppos\'e de composante neutre
r\'eductive ou de type (ssumult) pour $k$ non parfait; on note
$S=H^{\mult}$ le quotient maximal de $H$ qui est de type multiplicatif.
Soit $X=G/H$. Alors~:

\smallskip

(a) Il existe un isomorphisme $KD'(X) \simeq [\widehat T \to \widehat S]$
dans ${\cal D}(k)$ qui induit un isomorphisme fonctoriel
$$\Phi_X : \brr_{1,e} X \to H^1(k,[\widehat T \to \widehat S]) \, .$$

\smallskip

(b)  Notons
$$\ab^0=(\ab \circ \tau) : X(k) \to H^0_{\fp}(k,[S \to T])$$
la compos\'ee de $\tau$ avec l'application d'ab\'elianisation
$$\ab : H^0_{\fp}(k,[H \to G]) \to H^0_{\fp}(k,[S \to T]) \, .$$
Alors on a la formule~:
\begin{equation}\label{eq:compatibility}
\Phi_X(\alpha) \cup \ab^0(x)=-\alpha(x)
\end{equation}
pour tout $x \in X(k)$ et tout $\alpha \in \brr_{1,e} X$, o\`u
$\cup$ est l'accouplement (\ref{cupdef}).
\end{theo}

\rem Nous aurons notamment \`a appliquer le th\'eor\`eme sur un corps local de
caract\'eristique $p >0$, imparfait donc. L'hypoth\`ese faite sur
$H$ implique alors que $S=H^{\mult}$ est bien d\'efini.

\dem (a) Ceci a \'et\'e \'etabli par le premier auteur et
van Hamel \cite{borvh2} en caract\'eristique z\'ero. La m\^eme
preuve fonctionne en caract\'eristique quelconque une fois qu'on a
fait les hypoth\`eses suppl\'ementaires sur
$H$ du th\'eor\`eme~\ref{compatible}.

\smallskip

(b) La d\'emonstration de cette compatiblit\'e va occuper le reste de
cette section. Voir \'egalement l'appendice \ref{AppDem} pour une d\'emonstration
ind\'ependante du th\'eor\`eme~\ref{compatible} quand $H$ est de
type (ssumult).
On pourra \'egalement se reporter au lemme~4.5.1 de \cite{demthese}
pour une preuve de la compatibilit\'e analogue \`a celle que nous
donnons maintenant.

\smallskip

\begin{defi}
{\rm Une {\it paire de  $k$-groupes} est une paire $(G,H)$, o\`u $G$ est
un $k$-groupe lin\'eaire connexe lisse avec $\Pic \Gbar=0$
(suppos\'e r\'eductif si $k$ est non parfait), et
$H$ est un $k$-sous-groupe lisse de $G$ (non n\'ecessairement connexe, mais
de type (ssumult) ou avec $H^0$ r\'eductif pour $k$ non parfait).
}
\end{defi}

 Une paire $(G,H)$ d\'efinit  un  espace homog\`ene $X:=G/H$ avec un point marqu\'e $x^0=eH\in X(k)$,
 o\`u $e\in G(k)$ est  l'\'el\'ement neutre de $G$.
\bigskip

 Par un morphisme de paires $\phi\colon (G_1,H_1)\to (G_2,H_2)$
on entend  un  homomorphisme {\em surjectif} $\phi\colon G_1\to G_2$
 tel que $\phi(H_1)=H_2$.
 Si on pose $X_1=G_1/H_1$ et $X_2=G_2/H_2$, alors on a un   morphisme $\phi_*\colon (X_1, x^0_1)\to (X_2,x^0_2)$,
 o\`u $x^0_1$ et $x^0_2$ sont les points marqu\'es correspondants.

\bigskip

Soit   $(G,H)$ une paire de $k$-groupes    (avec $\Pic\Gbar=0$).
On choisit  un plongement $i\colon S:=H^{\mult}
\to Q$ dans un tore quasi-trivial $Q$.
On consid\`ere le plongement
$$
j\colon H \to G\times_k Q, \qquad h \mapsto (h,i(\mu(h))),$$
o\`u  $\mu\colon H \to S$ est  l'\'epimorphisme canonique.
On pose $G_Y=G\times_k Q$, $H_Y=i_*(H)$.
La paire $(G_Y,H_Y)$ d\'efinit un espace homog\`ene $Y=G_Y/H_Y=(G\times_k Q)/i_*(H)$ avec un point  marqu\'e  $y^0$.
 L'application de projection  $\pi\colon G_Y = G \times Q \to G$
est surjective et satisfait $\pi(H_Y) = H$, et  elle d\'efinit  donc un
 morphisme de paires  $\pi\colon (G_Y,H_Y)\to(G,H)$,
qui \`a son tour d\'efinit un morphisme de vari\'et\'es avec points marqu\'es  $\pi_*\colon (Y,y^0)\to (X,x^0)$.
 Notons encore $T_Y=G_Y^{\tor}$ et
$S_Y=H_Y^{\mult}$.
On remarque que l'application  $j_*\colon S_Y\to T_Y$ est injective et que $(Y,\pi_*)$ est un torseur sur $X$ sous $Q$.

On  construit  une nouvelle  paire $(G_Z,H_Z)$ comme suit :
 $G_Z=T_Y/j_*(S_Y)$, $H_Z=1$, $Z=G_Z$.
 Alors $Z$ est un $k$-tore, on d\'esigne  par $z^0$ son \'el\'ement neutre.
On a un morphisme de paires $\nu\colon (G_Y,H_Y)\to (G_Z, \{1\})$
et le morphisme induit  d'espaces homog\`enes
$\nu_*\colon (Y,y^0)\to (Z,z^0)$.
On a des quasi-isomorphismes
$$[\widehat{G_Y}\to\widehat{H_Y}]\to [\widehat{T_Y}\to \widehat{S_Y}]\to[\widehat{G_Z}\to 0].$$

On obtient des diagrammes
\begin{equation}\label{eq:main-diagram}
\xymatrix{
(G_Y, H_Y) \ar[d]^{\nu}\ar[r]^{\pi} & (G,H) & &(Y,y^0) \ar[d]^{\nu_*}\ar[r]^{\pi_*} & (X,x^0)\\
(G_Z,\{1\})                               &   & &(Z,z^0) \, .
}
\end{equation}

\begin{lem}\label{lem:Z}
Le th\'eor\`eme \ref{compatible} est valable pour l'espace principal homog\`ene $Z$ sous $G_Z$.
\end{lem}

\dem
Soit $Z$ un $k$-tore.
On a l'isomorphisme
$$
\Phi\colon\Brr_{1,e} Z\to H^2(k,Z)
$$
de \cite{borvh} et \cite{borvh2}.
On a aussi l'isomorphisme
$$
\Phi^S\colon\Brr_{1,e} Z\to H^2(k,Z)
$$
de Sansuc \cite{sansuc}, Lemme 6.9.
Ces isomorphismes diff\`erent par le signe, cf.  \cite{borvh2}, remarque 7.3 (sans preuve)
et l'appendice \ref{AppBor}, proposition \ref{prop:main-prop} du pr\'esent article.

On a un diagramme commutatif
$$
\xymatrix@C=0.7cm{
Z(k)\ar@{=}[d] &\times &\Br\, Z\ar[r]\ar[d]^{\Phi^S} &\Br\, k \ar@{=}[d]    \\
Z(k) &\times            &H^2(k,Z)\ar[r]                &\Br\, k,
}
$$
voir \cite{sansuc}, (8.11.2).
Ainsi pour $z\in Z(k)$, $\alpha\in\Br\, Z$ on a
$$
\alpha(z)=z\cup\Phi^S(\alpha)=-z\cup\Phi(\alpha),
$$
ce que prouve le  th\'eor\`eme \ref{compatible} pour $Z$.
\enddem

\begin{lem}\label{l:premain1}
Les diagrammes suivants sont commutatifs, toutes les fl\`eches 
verticales marqu\'ees  $(\cong)$ sont des isomorphismes,
et l'application $\pi_*\colon Y(k)\to X(k)$ est surjective :
{\small
\begin{equation}\label{eq:1}
\xymatrix@C=0.7cm{
X(k) \ar[r]^-{\ab^0} &H^0(k,X^\ab)  &  &H^0(k,X^\ab)
\ar@{}[r]|-{\times} &H^1(k,[\widehat{G}\to\widehat{H}]) \ar[d]^{\widehat{\pi}}  \ar[r]^-{\cup} &\Br k \ar@{=}[d]\\
Y(k) \ar[u]_{{\pi_*}} \ar[d]^{\nu_*} \ar[r]^-{\ab^0} &H^0(k,Y^\ab)\ar[u]_{{\pi}^\ab}\ar[d]_\cong^{\nu^\ab}
&  &H^0(k,Y^\ab) \ar[u]_{{\pi}^\ab}\ar[d]_\cong^{\nu^\ab}
                 \ar@{}[r]|-{\times} &H^1(k,[\widehat{G_Y}\to\widehat{H_Y}]) \ar[r]^-{\cup} &\Br k\\
Z(k) \ar[r]^-{\ab^0} &H^0(k,G_Z) &  &H^0(k,G_Z)
\ar@{}[r]|-{\times} &H^2(k,\widehat{G_Z})
     \ar[u]^\cong_{\widehat{{\nu}}}\ar[r]^-{\cup} &\Br k\ar@{=}[u] \, .
}
\end{equation}
}
\end{lem}

\dem
Le lemme r\'esulte de la fonctorialit\'e de $\ab^0$ et du cup produit et du fait
que le morphisme de complexes de $k$-groupes de type multiplicatif
$$
{\nu}^\ab\colon [ S_Y\to T_Y] \to [ 1\to G_Z]
$$
est un  quasi-isomorphisme.
L'application $\pi_*\colon Y(k)\to X(k)$ est surjective parce que $Y$ est un  torseur sur $X$ sous le  tore quasi-trivial $Q$.
\enddem

\begin{lem}\label{l:pre_main2}
Les diagrammes suivants sont commutatifs et 
toutes les fl\`eches verticales marqu\'ees  $(\cong)$ sont des isomorphismes :
{\small
\begin{equation}\label{eq:2}
\xymatrix@C=0.7cm{
X(k) \ar@{}[r]|-{\times} &\Brr_{1,m} X\ar[d]^{{{\pi}}^*} \ar[r]^{\langle\,,\,\rangle} &\Br k \ar@{=}[d]
& & &\Brr_{1,m} X\ar[d]^{\pi^*}\ar[r]_-\cong^-{\Phi_X} &H^1(k, [ \widehat{G}\to\widehat{H}])\ar[d]^{\widehat{\pi}} \\
Y(k) \ar[u]_{{\pi_*}} \ar[d]^{{\nu_*}} \ar@{}[r]|-{\times} &\Brr_{1,y_0} Y \ar[r]^{\langle\,,\,\rangle} &\Br k
& & &\Brr_{1,y_0} Y               \ar[r]_-\cong^-{\Phi_Y} &H^1(k, [\widehat{G_Y}\to\widehat{H_Y}])\\
Z(k)                \ar@{}[r]|-{\times} &\Brr_{1,w_0} Z \ar[u]^\cong_{{{\nu}}^*}\ar[r]^{\langle\,,\,\rangle} &\Br k\ar@{=}[u]
& & &\Brr_{1,w_0} W\ar[u]^\cong_{\nu^*}\ar[r]_-\cong^-{\Phi_Z} &H^2(k,  \widehat{G_Z})\ar[u]^\cong_{\widehat{\nu}} \, .
}
\end{equation}
}
\end{lem}

\dem
La commutativit\'e des diagrammes  r\'esulte de la fonctorialit\'e de $\langle\,,\,\rangle$ et de $\Phi$.
Comme $\widehat{\nu}$ est un isomorphisme, on voit que $\nu^*$ est un isomorphisme.
\enddem

\noindent
 {\bf Preuve du  Th\'eor\`eme \ref{compatible}.}
Par ``chasse au diagramme''  dans les diagrammes \eqref{eq:1} et \eqref{eq:2} on r\'eduit  le th\'eor\`eme au lemme \ref{lem:Z}.
\enddem

\section{Corps globaux et corps finis}

Rappelons tout d'abord un th\'eor\`eme d\'emontr\'e par Douai dans sa th\`ese
\cite{douai-these}~:

\begin{theo}[Douai] \label{douaitheo}
Soit $K$ un corps local. Soit $G$ un $K$-groupe semi-simple et soit $L$ un lien
localement repr\'esentable pour la topologie \'etale (resp. fppf)
par $G$. Alors toutes les \'el\'ements de $H^2(K,L)$ (resp. de $H^2_{\fp}(K,L)$)
sont neutres.
\end{theo}

Pour la preuve de ce r\'esultat,
on montre d'abord que le lien $L$ est globalement repr\'esentable par
une $K$-forme
de $G$ (\cite{douai}, Lemme~1.1), ce qui permet de se ramener au cas o\`u
$L$ est le lien canoniquement associ\'e \`a $G$. On utilise alors le fait que pour un
groupe semi-simple simplement
connexe $R$ de centre $Z$, l'application bord $H^1(K,R/Z) \to H^2(K,Z)$ est
surjective (\cite{douai}, Th. 1.1, ou \cite{boduke} Lemme 5.7 en
caract\'eristique z\'ero), ce qui permet de conclure via
\cite{gonz}, Cor. 3.8.

On en d\'eduit une g\'en\'eralisation d'un r\'esultat du
deuxi\`eme auteur (cf. \cite{demedin}, Prop. 2.18 et \cite{demthese}, Prop. 4.2.20
joints au fait que si $K$ est un corps local, alors $H^3_{\fp}(K,\mu)=0$
pour tout $K$-groupe fini $\mu$ \cite{adt}, Prop. 6.4)~:

\begin{prop} \label{abelocal}
Soit $K$ un corps local (de caract\'eristique quelconque). Soit $G$ un
$K$-groupe r\'eductif v\'erifiant $\pic \ov G=0$ et soit
$H$ un $K$-sous-groupe de $G$ v\'erifiant~: il existe des suites
exactes
$$1 \to H' \to H \to S \to 1$$
$$1 \to U \to H' \to H^{\ss} \to 1$$
avec~: $S$ de type multiplicatif, $H^{\ss}$ semi-simple,
et $U$ sous-groupe unipotent distingu\'e de $H$.
Soit $T=G^{\tor}$ le quotient torique maximal
de $G$.
Alors l'application canonique
$$\ab_K : H^0_{\fp}(K,[H \to G]) \to H^0_{\fp}(K,[S \to T])$$
est surjective.
\end{prop}

\rem Si $K$ est de caract\'eristique z\'ero, l'hypoth\`ese sur $H$
signifie exactement que $H$ est de type (ssumult). Noter aussi qu'on ne
demande pas ici que $S$ soit lisse.

\dem On utilise \cite{demedin}, lemme 2.1 et remarque 2.2, pour obtenir
un diagramme commutatif~:
$$
\xymatrix{
 H^0_{\fp}(K,G)\ar[r]\ar[d]^{f_1} &H^0_{\fp}(K,[H \to G])\ar[r]\ar[d]
 &H^1_{\fp}(K,H)\ar[r]\ar[d]^{f_2} &H^1_{\fp}(K,G)\ar[d]^{f_3}\\
H^0_{\fp}(K,T) \ar[r] &H^0_{\fp}(K,[S \to T])\ar[r]
&H^1_{\fp}(K,S)\ar[r] &H^1_{\fp}(K,T)
}
$$
o\`u la premi\`ere ligne est une suite exacte d'ensembles point\'es
et la deuxi\`eme une suite exacte de groupes ab\'eliens. Par ailleurs
on a une action du groupe $H^0_{\fp}(K,G)$ sur l'ensemble point\'e
$H^0_{\fp}(K,[H \to G])$, qui est compatible via le diagramme ci-dessus
avec l'action du groupe ab\'elien $H^0_{\fp}(K,T)$ sur
le groupe ab\'elien $H^0_{\fp}(K,[S \to T])$ par translation
(loc. cit., Prop. 2.6).
Ceci permet par "chasse au diagramme" de se ramener
\`a prouver les propri\'et\'es suivantes~: $f_1$ et $f_2$ sont surjectives,
$f_3$ a un noyau trivial. Comme $\pic \ov G=0$, le sous-groupe d\'eriv\'e
$G^{\ss}$ de $G$ est semi-simple simplement connexe, ce qui implique
$H^1(K,G^{\ss})=1$.
Ceci donne la surjectivit\'e de $f_1$ et la trivialit\'e du noyau de
$f_3$.

\smallskip

Il reste \`a prouver que $f_2$ est surjective. Soit $H_1=H/U$, on a
alors $S=H_1/H^{\ss}$ et il suffit de v\'erifier que les
fl\`eches $u_1 : H^1_{\fp}(K,H) \to H^1_{\fp}(K,H_1)$ et
$u_2 : H^1_{\fp}(K,H_1) \to H^1_{\fp}(K,S)$ sont surjectives.
Montrons que $u_2$ est surjective.
Soit $\xi\in H^1_\fp(K,S)$, alors l'obstruction \`a relever $\xi$ dans $H^1_\fp(K,H_1)$
est un \'el\'ement ${\rm ob}(\xi)\in H^2_\fp(K,L)$, o\`u $L$
est un lien localement repr\'esent\'e par $H^{\ss}$ (voir \cite{boduke}, Cor. 6.4 en caract\'eristique z\'ero et
la proposition IV.4.2.8(i) de \cite{giraud} dans le cas g\'en\'eral).
Par le  th\'eor\`eme~\ref{douaitheo} appliqu\'e \`a $H^{\ss}$, l'\'el\'ement ${\rm ob}(\xi)\in H^2_\fp(K,L)$ est neutre,
donc $\xi$ provient de $H^1_\fp(K,H_1)$ et la fl\`eche  $u_2$ est bien
surjective.

La surjectivit\'e de $u_1$ est quant \`a elle d\'emontr\'ee par Oesterl\'e :
voir \cite{oes}, IV.2.2, dans le cas o\`u $U$ et $H$ sont lisses.
 On v\'erifie que la m\^eme preuve fonctionne dans le cas plus g\'en\'eral
o\`u ne suppose plus que les groupes soient lisses,
en utilisant le fait que $H^2_{\fp}(K, \alpha_p)=0$. Cela conclut la preuve.

\enddem

\begin{theo} \label{globaltheo}
Soit $K$ un corps global de caract\'eristique $p \geq 0$.
Soit $G$ un $K$-groupe r\'eductif v\'erifiant $\pic \ov G=0$
et soit $H$ un $K$-sous-groupe lisse
de $G$ v\'erifiant les hypoth\`eses de la proposition~\ref{abelocal}.
On pose $X=G/H$ et on suppose que $X$ poss\`ede
une compactification lisse $X^c$. Soient $T=G^{\tor}$ le
quotient torique maximal
de $G$ et $S=H^{\mult}=H/H'$ le quotient maximal de $H$ qui est de type
multiplicatif.
Alors l'application $r$ de \S\,\ref{sec:4} induit un isomorphisme
$$\bra X^c/ \br K \to \Sha^1_{\omega}(K,[\widehat T \to \widehat S])$$
\end{theo}

\begin{rem}
Le th\'eor\`eme~\ref{theo transc} assure que, si de plus le
groupe $H$ est connexe (r\'eductif si $p > 0$), alors
$(\br X^c) \{p'\}=(\bra X^c) \{p'\}$. Sous
cette hypoth\`ese suppl\'ementaire, le th\'eor\`eme~\ref{globaltheo} donne
donc que
l'application $r$
induit un isomorphisme entre les sous-groupes de torsion premi\`ere
\`a $p$ de $\br X^c / \br K$ et $\Sha^1_{\omega}(K,[\widehat T \to \widehat S])$,
i.e. un isomorphisme
$$\left( \br X^c/ \br K \right)\{p'\} \to \Sha^1_{\omega}(K,[\widehat T \to \widehat S])\{p'\} \, .$$
En caract\'eristique nulle, si $H$ est connexe, on a donc un isomorphisme
$$\br X^c/ \br K \to \Sha^1_{\omega}(K,[\widehat T \to \widehat S]) \, .$$
\end{rem}

\bigskip

\paragraph{Preuve du th\'eor\`eme~\ref{globaltheo} : } Noter que
$H$ est de type (ssumult). On identifie
$\brr_1 X^c/\br K$ \`a $\brr_{1,e} \, X^c$. Soit $\alpha \in \brr_{1,e} \, X^c
\subset \brr_{1,e} \, X$. D'apr\`es le th\'eor\`eme~\ref{compatible}(a),
on a un isomorphisme $\Phi_X : \brr_{1,e} \, X \to
H^1(K,[\widehat T \to \widehat S])$, et la proposition~\ref{facilinclusion}
donne d\'ej\`a que l'image de $\brr_{1,e} \, X^c$ par
$\Phi_X$ contient $\Sha^1_{\omega}(K,[\widehat T \to \widehat S])$. Il reste
\`a montrer l'inclusion inverse.

\smallskip

Fixons $\alpha \in \brr_{1,e} X^c$.
Comme on l'a d\'ej\`a vu (preuve de la proposition~\ref{prop alteration}),
il existe un ensemble fini de places
$\Si$ de $K$ (contenant les \'eventuelles places archim\'ediennes)
tel que $\alpha(P_v)=0$ pour toute $v \not \in \Si$ et tout point local
$P_v \in X(k_v)$.

\smallskip

Soit $v \not \in \Si$. On sait que l'application
d'\'evaluation
$$\tau_v : X(K_v) \to H^0_{\fp}(K_v,[H \to G]); P_v \mapsto [[X]](P_v)$$
est bijective. Il en r\'esulte que
l'image de l'application d'\'evaluation
$$\ab^0 _v : X(K_v) \to H^0_{\fp}(K_v,[S \to T])$$
est aussi l'image de $$\ab_v : H^0_{\fp}(K_v,[H \to G])
\to H^0_{\fp}(K_v,[S \to T])$$
La proposition~\ref{abelocal} dit que $\ab_v$, et donc aussi
$\ab^0 _v$, est surjective.

\smallskip

Posons alors $a=\Phi_X(\alpha) \in H^1(K, [\widehat T \to \widehat S]))$.
D'apr\`es le th\'eor\`eme~\ref{compatible}(b), on obtient que pour
$v \not \in \Si$, la localisation
$a_v \in H^1(K_v,[\widehat T \to \widehat S]))$ de $a$ v\'erifie
$$a_v \cup \ab^0 _v(P_v)=\alpha(P_v)=0$$
pour tout $P_v \in X(K_v)$. Comme $\ab^0 _v$ est surjective, on obtient que
$a_v$ est
orthogonal \`a $H^0_{\fp}(K_v,[S \to T])$
pour l'accouplement local (\ref{cupdef}).
La proposition~\ref{perflocal} donne alors $a_v=0$; finalement on a bien
prouv\'e $\Phi_X(\alpha) \in \Sha^1_{\omega}(K,[\widehat T \to \widehat S])$
comme on voulait.

\enddem

\rem Si on ne  suppose pas l'existence d'une compactification lisse
pour $X$, une preuve similaire utilisant la proposition \ref{prop alteration} assure
que l'on a quand m\^eme pour $\ell \neq p$ un morphisme injectif entre les groupes de torsion $\ell$-primaires
$$(\brnr_1 X/ \br K)\{\ell \} \to \Sha^1_{\omega}(K,[\widehat T
\to \widehat S])) \{\ell \}$$
o\`u on a not\'e $\brnr_1 X:=\brnr X \cap \bra X$.

\begin{lem} \label{injh1}
Soit $K$ un corps. Soit $S \to T$ un morphisme de $K$-groupes de
type multiplicatif. On suppose que $T$ est un tore et on note $C^{\bullet}=
[\widehat T \to \widehat S]$ le complexe dual de $[S \to T]$.
Soit $K$ une extension de $k$ dans laquelle $k$ est
alg\'ebriquement ferm\'e. Alors la fl\`eche naturelle
$H^1(k, C^{\bullet}) \to H^1(K,C^{\bullet})$ est injective.
\end{lem}

\dem On consid\`ere le diagramme commutatif \`a lignes exactes
\begin{equation}\label{eq-ST}
\xymatrix{
H^1(k,\widehat T) \ar[r]\ar[d] &H^1(k, \widehat S) \ar[r]\ar[d] &H^1(k, C^{\bullet}) \ar[r]\ar[d] &H^2(k,\widehat T)\ar[d] \\
H^1(K,\widehat T) \ar[r] &H^1(K, \widehat S) \ar[r] &H^1(K, C^{\bullet}) \ar[r] &H^2(K,\widehat T) \, .
}
\end{equation}
On a une suite exacte
$$
1\to \Ga\to\Ga_K\to\Ga_k\to 1,
$$
o\`u $\Ga_K=\gal(\ov K/K)$, $\Ga=\gal(\ov K/K \kbar)$ et $\Ga_k=\gal(\kbar/k)=\gal(K \kbar/K)$.
Comme $T$ est un tore et l'action de $\Ga$ sur $\widehat{T}$ est triviale, on a  $H^1(\Ga,\widehat T)=\Hom(\Ga,\widehat{T})=0$.
On utilise alors les suites exactes de restriction-inflation
\def\That{{\widehat{T}}}
\def\Shat{{\widehat{S}}}
\begin{align*}
0\to &H^1(k,\Shat)\to H^1(K,\Shat)\to H^1(\Ga,\Shat)\\
0\to &H^1(k,\That)\to H^1(K,\That)\to H^1(\Ga,\That)=0\\
0\to &H^2(k,\That)\to H^2(K,\That)\to H^2(\Ga,\That) \, .
\end{align*}
On obtient alors que dans \eqref{eq-ST} la deuxi\`eme fl\`eche verticale est
injective,
la premi\`ere  est un isomorphisme et
la quatri\`eme est injective. Par chasse au diagramme, on obtient
le r\'esultat voulu.
\enddem

\begin{prop} \label{cycdeploye}
Soit $K$ un corps. Soit $S \to T$ un morphisme de $K$-groupes de
type multiplicatif. On suppose que $T$ est un tore et on note $C^{\bullet}=
[\widehat T \to \widehat S]$ le complexe dual de $[S \to T]$.
Soit $L$ une extension
finie galoisienne de $K$ qui d\'eploie $T$ et $S$, on pose
$\Ga=\gal(L/K)$.
On consid\`ere le sous-groupe $\Sha^1_{\omega,\alg}(\Gamma,C^{\bullet})\subset H^1(\Ga,C^{\bullet})$ constitu\'e des
\'el\'ements dont la restriction
\`a $H^1(\Ga_1, C^{\bullet})$ est nulle pour tout sous-groupe
cyclique $\Ga_1$ de $\Ga$.
 Alors~:
\smallskip

(a)l'homomorphisme d'inflation
$$
H^1(\Ga,C^{\bullet})\to H^1(K,C^{\bullet})
$$
est injectif.
 Si de plus $\Sha^2_{\omega,\alg}(\widehat T)=0$
(ex. $T$ quasi-trivial) ou si $S$ est un tore, alors cet homomorphisme
identifie $\Sha^1_{\omega,\alg}(\Gamma,C^{\bullet})$ avec
 $\Sha^1_{\omega,\alg}(K,C^{\bullet})$,
 et on a $\Sha^1_{\omega,\alg}(K,C^{\bullet})=
\Sha^1_{\omega}(K,C^{\bullet})$ si $K$ est un corps global.

\smallskip

(b) si $L$ est une extension cyclique de $K$, alors
$\Sha^1_{\omega,\alg}(K,C^{\bullet})=0$; si de plus $K$ est un corps
global, on a $\Sha^1_{\omega}(K,C^{\bullet})=0$.

\end{prop}

\rem Il ne semble pas clair que pour un complexe quelconque $[S \to T]$
de groupes de type multiplicatif sur un corps de nombres,
les groupes $\Sha^1_{\omega,\alg}(C^{\bullet})$ et
$\Sha^1_{\omega}(C^{\bullet})$ co\"{\i}ncident toujours; ce sera vrai
a posteriori (via le th\'eor\`eme~\ref{carzerotheo} ci-dessous)
pour un complexe provenant d'un espace homog\`ene comme dans le
th\'eor\`eme~\ref{globaltheo}.

\dem (a) On \'ecrit un diagramme commutatif \`a lignes exactes
{\small
$$
\xymatrix{
H^1(\Ga,\widehat T) \ar[r]\ar[d]^{\cong} &H^1(\Ga,\widehat S) \ar[r]\ar@{^{(}->}[d] &H^1(\Ga,[\widehat T \to\widehat S])\ar[r]\ar[d]
& H^2(\Ga,\widehat T) \ar[r]\ar@{^{(}->}[d] &H^2(\Ga,\widehat S)\ar[d]\\
H^1(K,\widehat T) \ar[r]\ar[d] &H^1(K,\widehat S)\ar[r]\ar[d] &H^1(K,[\widehat T \to\widehat S])\ar[r]\ar[d]
&H^2(K,\widehat T)\ar[r]\ar[d] &H^2(K,\widehat S)\ar[d]\\
0\ar[r] &H^1(L,\widehat S)\ar[r] &H^1(L,[\widehat T \to\widehat S])\ar[r]
 &H^2(L,\widehat T)\ar[r] &H^2(L,\widehat S) \, .
}
$$
}
Noter que $H^1(L,\widehat T)=0$ car $\Ga_L$ agit trivialement sur
le $\Z$-module libre de type fini $\widehat T$ (ce qui donne le z\'ero en
bas \`a gauche). La suite spectrale de Hochschild-Serre donne
les suites exactes de bas degr\'e
\begin{align*}
0\to H^1(\Ga,\That)&\to H^1(K,\That)\to H^1(L,\That)^{\Ga}\to H^2(\Ga,\That)\\
                                   &\to \ker[H^2(K,\That)\to H^2(L,\That)^{\Ga}]\to H^1(\Ga,H^1(L,\That)),\\
0\to H^1(\Ga,\Shat)&\to H^1(K,\Shat)\to H^1(L,\Shat)^{\Ga}\to H^2(\Ga,\Shat)\\
                                   &\to \ker[H^2(K,\Shat)\to H^2(L,\Shat)^{\Ga}]\to H^1(\Ga,H^1(L,\Shat)).
\end{align*}
On observe alors que la deuxi\`eme et la quatri\`eme colonnes sont exactes,
la fl\`eche $H^1(\Ga, \widehat T) \to H^1(K,\widehat T)$ est un isomorphisme, et
les fl\`eches $H^2(\Ga, \widehat T) \to H^2(K,\widehat T)$ et $H^1(\Ga, \widehat S) \to H^1(K,\widehat S)$ sont injectives.
Par chasse au diagramme,  on obtient alors que la fl\`eche
$H^1(\Ga,[\widehat T \to \widehat S]) \to H^1(K,[\widehat T \to \widehat S])$ est injective.

\smallskip

Comme l'action de $\Ga_L$ sur $\widehat T$ et $\widehat S$ est
triviale, on a
$$
\Sha^1_{\omega,\alg}(L,\widehat S)=0,\quad
\Sha^2_{\omega,\alg}(L,\widehat T)=0 \, .
$$
 La troisi\`eme ligne du diagramme
(et son analogue en rempla\c cant $L$ par $\Ga_1$, o\`u $\Ga_1$
est un sous-groupe
procyclique de $\Ga_L$) donne alors $\Sha^1_{\omega,\alg}(L,C^{\bullet})=0$.
Si $S$ est un tore, on sait de plus que la fl\`eche
$H^2(\Ga, \widehat S) \to H^2(K,\widehat S)$ est injective. Comme tout
\'el\'ement de $\Sha^1_{\omega,\alg}(K,C^{\bullet})$ est tu\'e dans
$H^1(L,C^{\bullet})$, une chasse au diagramme comme
ci-dessus montre alors qu'avec
l'une des hypoth\`eses suppl\'ementaires de (a), un tel \'el\'ement
provient de $H^1(\Ga, C^{\bullet})$. Ce qui pr\'ec\`ede montre alors
que $\Sha^1_{\omega,\alg}(\Ga,C^{\bullet})$
 s'identifie \`a $\Sha^1_{\omega,\alg}(K,C^{\bullet})$.
Si $K$ est un corps global, ceci donne
$\Sha^1_{\omega,\alg}(K,C^{\bullet})=
\Sha^1_{\omega}(K,C^{\bullet})$ par le th\'eor\`eme de \v Cebotarev.

\smallskip

(b) Si maintenant $L$ est une extension cyclique de $K$, alors on a en
particulier $\Sha^2_{\omega,\alg}(\widehat T)=0$ et le r\'esultat
d\'ecoule de la derni\`ere assertion de (a).
\enddem

\begin{theo} \label{finitheo}
Soit $k$ un corps fini de caract\'eristique $p$.
Soit $G$ un $k$-groupe r\'eductif. Soit $X$ un espace homog\`ene de
$G$, de stabilisateur g\'eom\'etrique $\ov H$.
on suppose que l'une des deux hypoth\`eses suivantes est satisfaite~:

\smallskip

(i) le $\kbar$-groupe $\ov H$ est connexe.

(ii) $\ov H$ est de type (ssumult) et $\pic \ov G=0$.

\smallskip

Alors si $X^c$ est une compactification lisse de $X$, on a $\bra X^c=0$.

Dans le cas (i), si $\ov H$ est r\'eductif, on a en outre $\br(X^c)\{p'\}=0$.
\end{theo}

\rem Dans le point (ii), l'hypoth\`ese $\pic \ov G = 0$ est n\'ecessaire,
comme on peut le voir dans l'exemple de la proposition \ref{prop fini ab},
o\`u le groupe $\ov H$ est fini de type multiplicatif et $\brnr_1(X)\{p'\} \neq 0$.

\dem Notons d\'ej\`a que $\br k=0$ (\cite{cogal}, section II.3.,
corollaire \`a la proposition 8).
D'autre part comme $G$ est lin\'eaire
connexe et $k$ fini, l'espace homog\`ene $X$ poss\`ede un point
rationnel par les th\'eor\`emes de Steinberg et de Springer
(\cite{cogal}, section III.2., corollaire au th\'eor\`eme 3,
ou \cite{springer}, 3.10).
Ainsi on peut supposer $X=G/H$, o\`u $H$ est un $k$-sous-groupe de
$G$ de type (ssumult) (on v\'erifie facilement que si $H \times_k \kbar$ est de
type (ssumult), alors $H$ est de type (ssumult)), qui est de plus suppos\'e
connexe si l'on ne fait pas l'hypoth\`ese $\pic \ov G=0$.

\smallskip

Soit alors $K=k(t)$. On commence par le cas o\`u
$\pic \ov G=0$.
Posons $T=G^{\tor}$. Comme $H$ est de type (ssumult),
on a une suite exacte de $k$-groupes lisses
$$ 1 \to H' \to H \to S \to 1$$
avec $S$ de type multiplicatif lisse, et $H'$ connexe lisse sans caract\`eres.
Comme $k$ est parfait, le groupe $H'$ est alors extension d'un $k$-groupe
semi-simple par un $k$-groupe unipotent $U$ connexe et lisse. Par ailleurs $U$ reste distingu\'e dans $H$ car $\ov U$
est invariant par tout $\kbar$-automorphisme de
$\ov H'$ (ce qui suffit vu que tous les groupes consid\'er\'es sont
lisses) en tant que radical unipotent de $\ov H'$.
Il en r\'esulte que le groupe $H_K :=H \times_k K$ satisfait
aux hypoth\`eses du th\'eor\`eme~\ref{globaltheo}.

\smallskip

Posons $C^{\bullet}:=[\widehat T \to \widehat S]$.
Alors d'apr\`es le th\'eor\`eme~\ref{compatible}(a), on a
$\bra X=H^1(k,C^{\bullet})$, ce qui donne une injection canonique
$\bra X^c\hookrightarrow H^1(k,C^{\bullet})$. De m\^eme on a une injection canonique
$\bra X^c_K /\br K \hookrightarrow H^1(K,C^{\bullet})$,
et le diagramme
$$
\xymatrix{
 \bra X^c\ar[r]\ar[d]  & H^1(k,C^{\bullet})\ar[d]\\
 \bra X^c_K/\br K\ar[r]  & H^1(K,C^{\bullet})
}
$$
est commutatif.
Ceci implique d'apr\`es le lemme~\ref{injh1} que la fl\`eche naturelle
$\bra X^c \to \bra X^c_K /\br K$ est injective.
Or $\bra X^c_K /\br K$
est isomorphe \`a $\Sha^1_{\omega}(K,C^{\bullet})$
d'apr\`es le th\'eor\`eme~\ref{globaltheo}. Comme toute extension finie
d'un corps fini est cyclique, la proposition~\ref{cycdeploye}(b)
s'applique et donne
$\Sha^1_{\omega}(K,C^{\bullet})=0$, d'o\`u $\bra X^c=0$.

\smallskip

Si on ne suppose plus $\pic \ov G=0$, on
observe (\cite {ctku2}, lemme 1.5)
que $X$ est isomorphe \`a un espace homog\`ene $G_1/H_1$,
o\`u $G_1$ est r\'eductif connexe quasi-trivial (en particulier
$\pic \ov G_1=0$) et $H_1$ est extension centrale de $H$ par un
$k$-tore $T_1$ ; en particulier $H_1$ reste lisse et connexe
et on est ramen\'e au cas pr\'ec\'edent.

\smallskip

Le dernier point de l'\'enonc\'e, sous l'hypoth\`ese $\ov H$ r\'eductif,
est alors une cons\'equence du th\'eor\`eme \ref{theo transc}.

\enddem

\rem L\`a encore, la preuve montre qu'on a $\brnr_1 X \{\ell \}=0$ pour
tout nombre premier $\ell \neq {\rm  Car} \, k$ sans supposer
l'existence d'une compactification lisse de $X$ : il suffit en effet de remplacer le th\'eor\`eme \ref{globaltheo}
par la remarque qui suit la preuve de ce m\^eme th\'eor\`eme \ref{globaltheo}.

Ce r\'esultat
est \'egalement valable sans supposer $G$ r\'eductif, en utilisant
le m\^eme argument que dans la preuve du th\'eor\`eme~\ref{carzerotheo}
ci-dessous (comme il s'agit d'un argument birationnel, il ne marche pas
pour la $p$-partie).

\section{Corps de caract\'eristique z\'ero}

Dans cette section, nous allons suivre la m\'ethode de
Colliot-Th\'el\`ene et Kunyavski\u\i\ \cite{ctku1}
pour \'etablir, par r\'eduction au cas des corps finis, le

\begin{theo} \label{carzerotheo}
Soit $k$ un corps de caract\'eristique z\'ero.
Soit $G$ un $k$-groupe lin\'eaire connexe
v\'erifiant $\pic \ov G=0$, on note
$T=G^{\tor}$ son quotient torique maximal.
Soit $X$ un espace homog\`ene de $G$ de stabilisateur g\'eom\'etrique $\ov H$,
on suppose que $\ov H$ est de type (ssumult) et on note
$S$ le $k$-groupe
de type multiplicatif canoniquement associ\'e \`a $\ov H$.
Soit $X^c$ une compactification lisse de $X$.
Alors on a une injection fonctorielle
$$\bra X^c/\br k \hookrightarrow
\Sha^1_{\omega,\alg}([\widehat T \to \widehat S])$$
qui est un isomorphisme si on suppose de plus $X(k) \neq \emptyset$
ou $H^3(k,\G)=0$.
\end{theo}

Rappelons que bien que $\ov H$ ne soit pas forc\'ement d\'efini sur $k$, il est
muni d'un $k$-lien qui induit une action de $\Ga_k$ sur le quotient
maximal de type multiplicatif $\ov S=\ov H^{\mult}$ de $\ov H$, et donc
une $k$-forme canonique $S$ de $\ov S$.

\dem
D'apr\`es le th\'eor\`eme~\ref{compatible}  (a)
(th\'eor\`eme 1 de \cite{borvh2}), on sait que
 $KD'(X)$ est isomorphe canoniquement \`a $[\widehat T \to \widehat S]$
dans la cat\'egorie d\'eriv\'ee,
donc il suffit de construire une injection
$$
\bra X^c/\br k \hookrightarrow\Sha^1_{\omega,\alg}(KD'(X)).
$$
On sait d\'ej\`a par la proposition~\ref{facilinclusion} qu'on a une
suite exacte (dont la derni\`ere fl\`eche est de plus z\'ero si
$X(k) \neq \emptyset$)
$$0 \to \brr_1 X/\br k \to H^1(k,KD'(X)) \to H^3(k,\G)$$
telle que de plus l'image de $\brr_1 X^c/\br k$ dans $ H^1(k,KD'(X))$ contienne
le noyau
$$\ker[\Sha^1_{\omega, \alg}(k,KD'(X))\to H^3(k,\G)],$$
voir
le diagramme \eqref{eq:diagBvH}
dans la preuve de la proposition \ref{facilinclusion}.
Comme cette image est contenue dans
celle de $H^1(k,KD'(X^c))=H^1(k,\pic \ov {X^c})$,
il ne reste plus
qu'\`a montrer que $$H^1(k,\pic \ov {X^c}) =
\Sha^1_{\omega,\alg}(k,\pic \ov {X^c}) .$$

\smallskip

Traitons d'abord le cas o\`u $G$ est r\'eductif et o\`u $X$
poss\`ede un $k$-point en d\'emontrant la

\begin{prop} \label{mainstep}
Avec les notations ci-dessus, supposons de
plus que le groupe $G$ est r\'eductif et que $X(k) \neq \emptyset$.
Alors pour tout
sous-groupe procyclique $\Ga_1$ de $\Ga_k$, on a
$$
H^1(\Ga_1, \pic \ov {X^c})=0.
$$
\end{prop}

\paragraph{Preuve de la proposition~\ref{mainstep} : } L'argument
va consister \`a se ramener au cas o\`u $k$ est fini; c'est
essentiellement le m\^eme que pour la proposition 2 de \cite{ctku1},
nous le reproduisons en entier pour le confort du lecteur.

\smallskip

On sait que $\ov G:= G \times_k \kbar$ est une $\kbar$-vari\'et\'e
rationnelle (\cite{borel}, corollaires 14.14 et 18.8).
De ce fait, $\ov X \simeq \ov G/\ov H$ (et donc aussi
$\ov X^c)$ est une vari\'et\'e
$\kbar$-unirationnelle puisqu'on a un $\kbar$-morphisme dominant
$\ov G \to \ov X$. Comme la caract\'eristique de $k$ est nulle,
ceci implique que pour tout $i >0$, les groupes $H^i(\ov X^c,\calo_{\ov X^c})$
(ou encore les groupes $H^i(X^c,\calo_{X^c})$)
sont nuls. Rappelons qu'on obtient
ce r\'esultat (bien connu) en se ramenant au cas $k=\CC$ par le principe de
Lefschetz, puis en utilisant le fait que $H^i(X^c,\calo_{X^c})=
H^0(X^c,\Omega_{X^c} ^i)$
via la th\'eorie de Hodge (c'est cet argument qui ne marche pas en
caract\'eristique $p$). La nullit\'e de ces groupes est alors assur\'ee
par le fait que $\ov X^c$ est domin\'e par un espace projectif, car
on a $H^i({\bf P}^r,\calo_{{\bf P}^r})=0$ pour tous $i,r >0$ via
\cite{hart}, th\'eor\`eme III.2.7 et th\'eor\`eme III.5.1(b),(d).
La propri\'et\'e
que $\ov X^c$ est unirationnelle implique \'egalement que $\pic \ov X^c$ est
sans torsion (par exemple parce que d'apr\`es un th\'eor\`eme de
Serre \cite{serreunir},
le groupe fondamental g\'eom\'etrique $\pi_1(\ov X^c)$ est nul, ce qui donne
que le groupe $(\pic \ov X^c)[n] \simeq H^1(\ov X^c,\Z/n)$ est nul); le
fait que $H^1(X^c,\calo_{X^c})=0$ implique que $\pic \ov X^c$ est
de type fini (car \'egal au groupe de N\'eron-Severi de $\ov X^c$).

\smallskip

\def\ktil{{\tilde{k}}}
\def\Xck{{{\vphantom{\ov X} X}^c_\ktil}}

Soit maintenant $\tilde k/k$
une extension finie galoisienne d\'eployant $\pic \ov X^c$,
posons $\Ga=\gal(\ktil/k)$. Le
fait que le groupe ab\'elien $\pic \ov X^c$ soit libre et de type fini implique
$H^1(\ktil, \pic \ov X^c)=0$.
Il suffit alors pour montrer la proposition de v\'erifier que pour tout
sous-groupe cyclique $\gamma \subset \Ga$, on a $H^1(\gamma,\pic \Xck)=0$
(o\`u $\Xck=X^c \times_k \ktil$). On \'ecrit
$\ktil=k[t]/P(t)$, o\`u $P \in k[t]$
est un polyn\^ome
irr\'eductible. Comme $X \hookrightarrow X^c$ est un morphisme de $k$-sch\'emas de
type fini,
on peut trouver un corps $k_0 \subset k$ de type fini sur $\Q$ sur lequel
$G,X,X^c$, et le plongement $X \to X^c$ sont d\'efinis; comme
$X(k) \neq \emptyset$ par hypoth\`ese, on peut \'egalement (quitte \`a remplacer
$k_0$ par le corps r\'esiduel de $x \in X_{k_0}$ d\'efini comme l'image d'un
$k_0$-morphisme $\spec k \to X_{k_0}$)
supposer que
$X(k_0) \neq \emptyset$ (ce qui implique que $X_{k_0}$ est isomorphe \`a
un quotient $G_{k_0}/H_{k_0}$, o\`u $H_{k_0}$ est un $k_0$-groupe de
type (ssumult)).
Quitte \`a agrandir $k_0$ dans $k$,
on peut aussi supposer que $P(t) \in k_0[t]$ et que $\tilde k_0:=k_0[t]/P(t)$
d\'eploie $\pic \ov X^c$ (qui est de type fini).
On peut enfin supposer
que $\tilde k_0$ est
galoisienne sur $k_0$~: en effet si $a \in \tilde k_0$ est une racine de $P$,
on a $\tilde k_0=k_0[a]$ et $\tilde k=k[a]$; toute racine $a_i$ de $P$ s'\'ecrit
alors $Q_i(a)$ avec $Q_i \in k[t]$ (parce que $\tilde k$
est galoisienne sur $k$),
et il suffit alors de prendre $k_0$ assez grand pour que tous les
polyn\^omes $Q_i$ soient dans $k_0[t]$. L'irr\'eductibilit\'e de $P$
sur $k_0$ implique alors que $\gal(\tilde k/k)$ est canoniquement isomorphe
\`a $\gal(\tilde k_0/k_0)$, et par ailleurs
les inclusions
$$
\pic X^c_{\tilde k_0} \to \pic \Xck \to \pic \ov X^c
$$
sont des isomorphismes. Ainsi le $\gal(\tilde k/k)$-module
$\pic \Xck$ est
canoniquement isomorphe au $\gal(\tilde k_0/k_0)$-module
$\pic {X}^c_{\tilde k_0}$.
Pour simplifier les notations, on notera dans la suite
$k=k_0$ et $\tilde k=\tilde k_0$ (en particulier
on s'est ramen\'e au cas o\`u $X=G/H$ avec
$H$ de type (ssumult)).

\smallskip

On prend alors un mod\`ele entier de toute la situation~: comme $k$
est
un corps de type fini sur $\Q$, on peut trouver un anneau int\`egre et
r\'egulier $A$ de type fini sur $\Z$, de corps des fractions $k$, tel
que la fermeture int\'egrale $\til A$
de $A$ dans $\tilde k$ soit finie, \'etale, galoisienne
de groupe $\Ga$ sur $A$ , et des $A$-sch\'emas de type fini
${\cal G}$, ${\cal H}$,
${\cal X}={\cal G}/{\cal H}$, ${\cal X}^c$
de fibres g\'en\'eriques respectives $G$, $H$, $X$, $X^c$, tels que~:

(i) ${\cal G}$ est un $A$-sch\'ema en groupes r\'eductif (en particulier
lisse \`a fibres connexes).

(ii) ${\cal H}$ est un $A$-sch\'ema en groupes lisse et ${\cal X}$ est
un $A$-sch\'ema lisse.

(iii) ${\cal X}^c$ est propre et lisse sur $A$ et on a une $A$-immersion
ouverte
${\cal X} \to {\cal X}^c$ qui \'etend $X \to X^c$.

(iv) Pour tout point $x$ de $\spec A$, de corps r\'esiduel $\kappa(x)$, la fibre
${\cal X}_x \subset {\cal X}^c_x$ est une compactification lisse
de ${\cal X}_x$ sur $\kappa(x)$.

\smallskip

On utilise alors que $H^1(X^c,\calo_{X^c})=H^2(X^c,\calo_{X^c})=0$. Le
th\'eor\`eme de semi-continuit\'e (\cite{hart}, III.12.8)
nous permet de supposer (quitte \`a
restreindre encore $\spec A$) que pour tout point $x$ de $\spec A$, on a encore
$H^1({\cal X}^c_x,\calo_{{\cal X}^c_x})=H^2({\cal X}^c_x,
\calo_{{\cal X}^c_x})=0$.
Les th\'eor\`emes de Grothendieck sur les faisceaux inversibles
(\cite{fga}, paragraphe 5; voir aussi \cite{dhthese}, Prop. 3.4.2)
disent
alors qu'on a un isomorphisme de sp\'ecialisation $\pic X^c \simeq
\pic {\cal X}^c_x$; plus pr\'ecis\'ement le foncteur
$\underline{\pic}_{{\cal X}^c/A}$ est repr\'esent\'e par un $A$-groupe
constant tordu
localement libre pour la topologie \'etale, et on peut supposer (quitte
\`a restreindre encore $A$ et \`a augmenter $\tilde k$ et $\til A$)
que ce $A$-groupe
est d\'eploy\'e (autrement dit devient isomorphe \`a $\Z^n$) sur $\spec \til A$.

Quitte \`a restreindre encore $\spec A$, on peut aussi supposer
que les fibres de ${\cal G}$ ont un groupe de Picard g\'eom\'etrique nul
(en effet sur un ouvert $U$ de $\spec A$ on a un rev\^etement simplement
connexe ${\cal G}^{\sc}$ de ${\cal G}^{\ses}$, qui est un rev\^etement fini
\'etale de ${\cal G}^{\ses}$ trivial au point g\'en\'erique, donc sur $U$ tout
entier). De m\^eme on peut supposer que les fibres g\'eom\'etriques de
${\cal H}$ sont de type (ssumult) (car $H$ est extension d'un $k$-groupe
de type multiplicatif par un groupe unipotent et ces propri\'et\'es
s'\'etendent sur un certain ouvert de $\spec A$).

\smallskip

Soit maintenant $\gamma$ un sous-groupe cyclique de $\Ga$ (correspondant \`a
une sous-extension $k_1 \supset k$ de $\tilde k$) et $A_1:=\til{A}^{\gamma}$
le sous-anneau de $\til A$ des invariants pour l'action de $\gamma$.
Le morphisme $\rho : \spec \til A \to \spec A_1$ est un rev\^etement \'etale cyclique
de groupe $\gamma$, et
les anneaux $\til A$ et $A_1$ sont int\`egres, r\'eguliers, de type fini sur $\Z$.
La version de Serre du th\'eor\`eme de \v Cebotarev (\cite{ceboserre}, th\'eor\`eme 7)
dit alors qu'il existe un
id\'eal maximal (en fait une infinit\'e) $m$ de $A_1$ tel que la fibre en $m$
de $\rho$ soit du type $\spec \til F$, o\`u $\til F$
est un corps fini (extension finie
du corps $F_1:=A_1/m$), autrement dit il y a un seul point $\til m$ de $\spec \til A$
au-dessus de $m$ et on a
$\gamma:=\gal(\tilde k/k_1)=\gal(\til F/F_1)$.
Soit $Y={\cal X}^c \times_A F_1$ la fibre de
${\cal X}^c \times_A A_1$ en $m$.
Comme $\underline{\pic}_{{\cal X}^c/A} \times_A \til A$
est un sch\'ema en groupes constant, sa fibre g\'en\'erique est isomorphe \`a
sa fibre en $\til m$ et on obtient
$$
H^1(\gamma,\pic \Xck)=H^1(\gamma, \pic (Y \times_F \til F)).
$$

Or $\gal(\ov F/\til F)$ agit trivialement sur $\pic \ov Y$, donc
$H^1(\gamma, \pic (Y \times_F \til F))$
est un sous-groupe de $H^1(F_1,\pic \ov Y)$. Comme $Y$ est
une compactification lisse d'un espace homog\`ene d'un $F_1$-groupe
lin\'eaire ${\cal G}_m$
connexe de groupe de Picard g\'eom\'etrique nul et \`a stabilisateur ${\cal H}_m$
de type (ssumult), le th\'eor\`eme \ref{finitheo} dit que $H^1(F_1,\pic \ov Y)=\Brr_1 Y=0$,
donc $H^1(\gamma,\pic \Xck)=0$,
ce qui ach\`eve la preuve de la proposition. 
\enddem

\def\isoto{\overset{\sim}{\to}}
\def\Xp{{\vphantom{\ov X} X}}
\paragraph{Fin de la preuve du th\'eor\`eme~\ref{carzerotheo} : }

Supposons d'abord que $G$ est r\'eductif (mais pas forc\'ement que
$X(k) \neq \emptyset$). Pour se ramener au cas $X(k) \neq \emptyset$, on utilise
l'``astuce du point g\'en\'erique'' (comme dans \cite{ctku1}, p. 6).
Soit $K$ le corps des fonctions de $X$, alors $X(K) \neq \emptyset$.
Soit $\ov K$ une cl\^oture alg\'ebrique de $K$,
alors $\kbar$ s'identifie avec la fermeture alg\'ebrique de $k$ dans $\ov K$,
et on a $\kbar\cap K=k$ (les corps $K$ et $\kbar$ sont lin\'eairement disjoints
sur $k$ car $X$ est g\'eom\'etriquement int\`egre).
 Soit comme ci-dessus
$\tilde k$ une extension galoisienne finie de $k$ dans $\kbar$
d\'eployant $\pic  \Xp^c_\kbar$,
alors $\gal(\tilde k/k)$ agit sur $\pic  \Xp^c_\kbar$.
On pose $\til K=K \tilde k\subset\ov K$, alors $\til K$ d\'eploie $\pic  \Xp^c_{\ov K}$, et
$\gal(\til K/K)$ agit sur $\pic  \Xp^c_{\ov K}$.
On a un isomorphisme canonique
$\varepsilon\colon \gal(\til K/K)\isoto \gal(\tilde k/k)$
et un isomorphisme canonique $\varepsilon$-\'equivariant
$\pic \Xp^c_\kbar \isoto \pic \Xp^c_{\ov K}$.
Comme $X(K)\ne\emptyset$, par  la proposition~\ref{mainstep},
pour chaque sous-groupe cyclique $\Gamma_1\subset\gal(\til K/K)$ on a
$H^1(\Gamma_1,\pic \Xp^c_{\ov K})=0$.
On conclut que pour chaque sous-groupe cyclique
$\gamma_1\subset\gal(\tilde k/k)$ on a
$H^1(\gamma_1,\pic \Xp^c_\kbar)=0$.
Ceci \'etablit le th\'eor\`eme~\ref{carzerotheo} pour $G$ r\'eductif.

Passons enfin au cas g\'en\'eral o\`u on ne suppose plus $G$ r\'eductif.
Soient $G^u$ le radical unipotent de
$G$ et $Y=X/G^u$ (voir \cite{bocrelle}, lemme 3.1 pour la repr\'esentabilit\'e de ce quotient).
Alors $Y$ est un espace homog\`ene sous le $k$-groupe
r\'eductif $G':=G/G^u$, de stabilisateur g\'eom\'etrique $\ov H'=\ov H/ (\ov H
\cap \ov G_u)$, pour lequel $(\ov H')^{\mult}=\ov H^{\mult}$ vu que
$(\ov H \cap \ov G_u)$ est un groupe unipotent.
On a encore $\pic \ov G'=0$ via \cite{sansuc} corollaire 6.11,
et $(G')^{\tor}=G^{\tor}$ par d\'efinition de $G^u$.
Il suffit donc, pour se ramener
au cas $G$ r\'eductif, de
montrer que $\brnr_1 Y=\brnr_1 X$. Or ceci r\'esulte de ce que
$X$ est stablement $k$-birationnel \`a $Y$, par le m\^eme argument
qu'\`a la fin de la deuxi\`eme \'etape dans la preuve du
th\'eor\`eme~\ref{theo transc}.

\enddem

\rems -Quand $H$ est connexe, la m\^eme m\'ethode donne encore
$\bra X^c/\br k \simeq  \Sha^1_{\omega,\alg}(KD'(X))$
sans supposer que $\pic \ov G=0$, mais $KD'(X)$ n'a plus en g\'en\'eral
une description aussi agr\'eable que quand $\pic \ov G=0$, voir \cite{DemBrauer}, th. 0.1.

\smallskip

-Si on ne suppose plus $k$ de caract\'eristique z\'ero, mais si on
suppose de plus $G$ r\'eductif (et toujours $\ov H$ de type
(ssumult)), la m\^eme m\'ethode fonctionne {\it \`a
condition de savoir qu'on a
$H^1(X^c,\calo_{X^c})=0=H^2(X^c,\calo_{X^c})=0$}. Or ceci est une
question ouverte en caract\'eristique $p$ pour les vari\'et\'es
s\'eparablement unirationnelles comme $X^c$.

\smallskip

On d\'eduit du th\'eor\`eme \ref{carzerotheo} le corollaire suivant, via le th\'eor\`eme \ref{theo transc} :
\begin{cor}
Avec les notations et les hypoth\`eses du th\'eor\`eme \ref{carzerotheo}, si on suppose en outre que $\ov H$ est connexe,
 alors on a une injection fonctorielle
$$\br X^c/\br k \hookrightarrow \Sha^1_{\omega,\alg}([\widehat T \to \widehat S])$$
qui est un isomorphisme si on suppose de plus $X(k) \neq \emptyset$
ou $H^3(k,\G)=0$.
\end{cor}

\smallskip

Le th\'eor\`eme~\ref{carzerotheo} n'est plus valable si on ne suppose
pas $\ov H$ de type (ssumult). On trouve dans \cite{demannalen} un contre-exemple avec
$G={\rm SL}_n$ et $\ov H$ fini non commutatif. Nous allons pour conclure
en d\'eduire un contre-exemple similaire avec $\ov H$ extension d'un groupe
fini commutatif par un tore, ce qui montre que l'hypoth\`ese
que le groupe des composantes connexes $\ov H/\ov H ^0$ est commutatif
n'est pas non plus suffisante.

\smallskip

Soit $p$ un nombre premier. On consid\`ere le groupe fini $E$
consid\'er\'e dans la proposition~6.1 de \cite{demannalen}. Il est
donn\'e par la pr\'esentation
$$E=\langle x,y,z : x^{p^2}=y^{p^2}=z^{p^2}=1; [x,y]=z^p \rangle \, .$$
C'est une
extension centrale de $F=E^{\ab}=
\Z/p^2 \times \Z/p^2 \times \Z/p$ par $A=D(E)=\Z/p$.
Soit
$k=\Q(\zeta_p)$, on regarde $E$ comme un groupe alg\'ebrique fini
(constant) sur $k$.
On plonge $\Z/p=\mu_p$ dans $\G$, et on d\'efinit un
$k$-groupe alg\'ebrique $H$ par $H:=(E \times \G)/j(A)$, o\`u
$j : A \to E \times \G$ est le plongement diagonal (noter que $j(A)$ est
bien distingu\'e dans $E \times \G$ car $A$ est central dans $E$). On plonge
ensuite $H$ dans un groupe ${\rm SL}_n$.

\begin{prop} \label{contrexprop}
Soit $X$ l'espace homog\`ene d\'efini par
$X={\rm SL}_n/H$. Alors $\brr_{1,e} X=
H^1(k,\widehat H^{\mult})$, mais $\brnr_{1,e} X$ contient strictement
$\Sha^1_{\omega}(\widehat H^{\mult})$.
\end{prop}

\dem Le fait que $\brr_{1,e} X=
H^1(k,\widehat H^{\mult})$ r\'esulte du th\'eor\`eme 1
de \cite{borvh2} joint au fait que ${\rm SL}_n$ est sans caract\`eres
et toute fonction inversible sur $X$ est constante.

\smallskip

On a
un diagramme commutatif, dont les deux premi\`eres lignes sont
exactes et correspondent \`a des extensions centrales~:
$$
\xymatrix{
1 \ar[r] &A \ar[r]\ar[d] &E \ar[r]\ar[d] &F \ar[r]\ar@{=}[d] &1 \\
1 \ar[r] &\G \ar[r] &H \ar[r] &F \ar[r] &1 \, .
}
$$
Comme le groupe
$F$ est de type multiplicatif (il est fini commutatif et $k$ est de
caract\'eristisque $0$), la fl\`eche $H \to F$ du diagramme se
factorise par une fl\`eche $H^{\mult} \to F$ d'o\`u une suite exacte
$$ \G \to H^{\mult} \to F \to 1 \, .$$
Il en r\'esulte que la fl\`eche induite $H^{\mult} \to F$ a pour noyau
un tore $T$ isomorphe \`a $\G$ puisque ce noyau est un quotient non trivial
de $\G$ ($H^{\mult}$ ne peut pas \^etre fini car $H/E=\G/A$ est isomorphe
\`a $\G$).

\smallskip

Ainsi $H^{\mult}$ est isomorphe \`a $T \times F$
(il n'y a pas de groupe de type multiplicatif extension non triviale d'un
groupe de type multiplicatif par un tore d\'eploy\'e), ce qui implique
que pour toute extension de corps $K$ de $k$, les fl\`eches canoniques
$$H^1(K, H^{\mult}) \to H^1(K,F) ; \quad H^1(K,\widehat F) \to H^1(K, \widehat
H^{\mult})$$
sont des isomorphismes. De ce fait, pour tout compl\'et\'e $k_v$ de
$k$, la fl\`eche
$$\im [H^1(k_v,H) \to H^1(k_v,H^{\mult})] \to
\im [H^1(k_v,H) \to H^1(k_v,F)]$$
est un isomorphisme. D'apr\`es \cite{demannalen} (preuve de la
proposition~6.1), il existe $a \in H^1(k,\widehat F)$
dont, pour une infinit\'e de $v$,
la localisation $a_v \in H^1(k_v,\widehat F)$ est non nulle, mais tel que
pour presque toute place $v$ on a $a_v$
orthogonale
\`a $I_v:=\im [H^1(k_v,E) \to H^1(k_v,F)]$, laquelle contient
$\im [H^1(k_v,H) \to H^1(k_v,F)]$ car la fl\`eche $H^2(k_v,A)=(\br k_v)[p]
\to H^2(k_v,\G)=\br k_v$ est injective. Soit alors $b$ l'image de
$a$ dans $H^1(k,\widehat H^{\mult})$; on a, pour une infinit\'e de
places $v$, que la localisation $b_v \in H^1(k_v,\widehat H^{\mult})$
est non nulle, et aussi que $b_v$ est orthogonale \`a
$\im [H^1(k_v,H) \to H^1(k_v,H^{\mult})]$ pour presque toute place $v$.
En particulier $b$ n'est pas dans $\Sha^1_{\omega}(\widehat H^{\mult})$.

\smallskip

D'apr\`es la formule de compatibilit\'e (voir th\'eor\`eme~\ref{compatible}(b))
et le fait que la fl\`eche d'\'evaluation $X(k_v) \to H^1(k_v,H)$
(associ\'ee au $H$-torseur ${\rm SL}_n \to X$) soit surjective
(ceci r\'esulte de ce que $H^1(k_v,{\rm SL}_n)=0$),
on obtient alors que pour presque toute
place $v$, on a $b(P_v)=0$ pour tout $k_v$-point $P_v$ de $X$, ce
qui montre que $b \in \brnr X$ via le th\'eor\`eme~2.2.1 de \cite{dhthese}.

\enddem

On peut \'egalement modifier l\'eg\'erement la construction pr\'ec\'edente pour v\'erifier que
dans le cas d'un corps fini, l'hypoth\`ese $\pic \ov G = 0$ dans le th\'eor\`eme \ref{finitheo} (ii)
est n\'ecessaire (voir la remarque qui suit le th\'eor\`eme \ref{finitheo}).

\begin{prop} \label{prop fini ab}
Il existe un corps fini $\F$ de caract\'eristique $p \geq 3$, un groupe semi-simple connexe $G$
(de groupe fondamental $\mu_2$) sur $\F$ et un sous-$\F$-groupe $H$ de $G$,
tel que $H$ soit un groupe fini constant commutatif d'ordre $16$, isomorphe \`a $(\Z / 2 \Z)^2 \times \Z / 4 \Z$,
de sorte que si $X = G / H$, on a
$$\brnr_1(X)\{2\} \neq 0 \, .$$
\end{prop}

\dem
On remarque d'abord que dans la proposition 6.1 de \cite{demannalen}, dans le cas d'un $2$-groupe,
on peut remplacer le groupe le groupe $E$ d'ordre $64$ par le groupe fini $H_0$ d'ordre $32$ dont une pr\'esentation est
$$H_0 =\langle x,y,z : x^2 = y^4 =z^4 = 1; [x,y]=z^2 \rangle \, .$$
Il s'agit par exemple du groupe d'ordre $32$ de num\'ero 24 dans la classification de GAP.

La preuve de la propostion 6.1 dans \cite{demannalen} assure que pour toute $\Q$-repr\'esentation fid\`ele
$\rho : H_0 \to \SL_{n,\Q}$, le groupe $\brnr_1(X_{\rho}) \{2 \}$ est strictement plus gros
que le sous-groupe de $\br_1(X_{\rho})\{2\}$ form\'e des \'el\'ements localement constants en presque tout place,
avec $X_{\rho} := \SL_{n, \Q} / \rho(H_0)$.

On consid\`ere la $\Q$-repr\'esentation $\rho_0$ suivante : on d\'efinit $\rho_0 : H_0 \to \SL_{8, \Q}$ par les formules
$$\rho_0(x) :=
\left(\begin{array}{cccccccc}
0 & 0 & 1 & 0 & 0 & 0 & 0 & 0 \\
0 & 0 & 0 & 1 & 0 & 0 & 0 & 0 \\
1 & 0 & 0 & 0 & 0 & 0 & 0 & 0 \\
0 & 1 & 0 & 0 & 0 & 0 & 0 & 0 \\
0 & 0 & 0 & 0 & 0 & 0 & 1 & 0 \\
0 & 0 & 0 & 0 & 0 & 0 & 0 & 1 \\
0 & 0 & 0 & 0 & 1 & 0 & 0 & 0 \\
0 & 0 & 0 & 0 & 0 & 1 & 0 & 0
\end{array}\right), \textup{  }
\rho_0(y) :=
\left(\begin{array}{cccccccc}
0 & -1 & 0 & 0 & 0 & 0 & 0 & 0 \\
1 & 0 & 0 & 0 & 0 & 0 & 0 & 0 \\
0 & 0 & 0 & 1 & 0 & 0 & 0 & 0 \\
0 & 0 & -1 & 0 & 0 & 0 & 0 & 0 \\
0 & 0 & 0 & 0 & 1 & 0 & 0 & 0 \\
0 & 0 & 0 & 0 & 0 & 1 & 0 & 0 \\
0 & 0 & 0 & 0 & 0 & 0 & -1 & 0 \\
0 & 0 & 0 & 0 & 0 & 0 & 0 & -1
\end{array}\right)
$$

$$\textup{et } \rho_0(z) :=
\left(\begin{array}{cccccccc}
0 & -1 & 0 & 0 & 0 & 0 & 0 & 0 \\
1 & 0 & 0 & 0 & 0 & 0 & 0 & 0 \\
0 & 0 & 0 & -1 & 0 & 0 & 0 & 0 \\
0 & 0 & 1 & 0 & 0 & 0 & 0 & 0 \\
0 & 0 & 0 & 0 & 0 & -1 & 0 & 0 \\
0 & 0 & 0 & 0 & 1 & 0 & 0 & 0 \\
0 & 0 & 0 & 0 & 0 & 0 & 0 & -1 \\
0 & 0 & 0 & 0 & 0 & 0 & 1 & 0
\end{array}\right) \, .
$$

On trouvera dans \cite{chkp} une repr\'esentation complexe de $H_0$,
tr\`es similaire \`a la repr\'esentation $\rho_0$ (dans \cite{chkp}, le groupe $H_0$ s'appelle (16;24)).

On v\'erifie facilement que $\rho_0$ d\'efinit une repr\'esentation fid\`ele de $H_0$ dans $\SL_{8, \Q}$.

Or $\rho_0(z^2) = -I_8$ est central dans $\SL_{8,\Q}$, donc on dispose du diagramme commutatif suivant,
o\`u $G := \SL_{8, \Q} / \{ \pm I_8 \}$ et $H := H_0^{\textup{ab}} = H/\langle z^2 \rangle$ :
\begin{displaymath}
\xymatrix{
1 \ar[r] & \langle z^2 \rangle \ar[r] \ar[d]^= & H_0 \ar[r] \ar[d]^{\rho_0} & H \ar[r] \ar[d] & 1 \\
1 \ar[r] & \mu_{2, \Q} \ar[r] & \SL_{8, \Q} \ar[r] & G \ar[r] & 1 \, .
}
\end{displaymath}
Ce diagramme permet d'identifier $X_{\rho_0} = \SL_{8,\Q}/H_0$ avec le quotient $G / H$ (comme $\Q$-vari\'et\'es).

Enfin, on a bien $H \cong (\Z / 2 \Z)^2 \times \Z / 4 \Z$, et $G$ est lin\'eaire connexe semi-simple, de groupe fondamental $\mu_2$.

On \'etend ensuite $H$, $G$, $\rho_0$ et $X_{\rho_0}$ sur un ouvert suffisamment petit de $\spec(\Z)$, de bonne r\'eduction.
Supposons que pour tout point ferm\'e $s$ de cet ouvert, le groupe de Brauer non ramifi\'e de la fibre de $X_{\rho_0}$ au-dessus
de $s$ (qui est un espace homog\`ene d'un groupe semi-simple connexe \`a stabilisateur $H$ sur le corps fini $k(s)$) soit trivial.
Alors la preuve de la proposition 2 de \cite{ctku1} assure que $\brnr_1(X_{\rho_0})$ coincide avec le sous-ensemble de $\br_1(X_{\rho_0})$
form\'es des \'el\'ements localement constants en presque toute place. Cela contredit la proposition 6.1 de \cite{demannalen}.

Donc il existe bien un corps fini $\F$ et un espace homog\`ene sur $\F$ de la forme $X = G/H$, avec $G$ semi-simple (connexe)
et $H$ ab\'elien fini constant d'ordre $16$ tels que $\brnr_1(X)\{2\} \neq 0$.
\enddem

\rem La preuve assure que l'on peut aussi \'ecrire $X$ sous la forme $X = \SL_{8, \F} / H_0$,
o\`u $H_0$ est un groupe fini constant d'ordre $32$.
Cela fournit en particulier un exemple de groupe fini $H_0$ d'ordre $32$ tel que le corps des fonctions de $X$
ne soit pas une extension transcendante pure du corps fini $\F$.
Cela donne en particulier un contre-exemple au probl\`eme de Noether sur un corps fini.

\appendix
\renewcommand{\thetheo} {\Alph{section}.\arabic{theo}}

\section{Compatibilit\'e pour un tore}\label{AppBor}

\def\muB{{\mu_{\rm BvH}}}
\def\lambdaS{{\lambda_{\rm S}}}

\def\deltaB{{\delta_{\rm BvH}}}
\def\deltaN{{\delta_{\rm nat}}}

\def\That{{\widehat{T}}}
\def\Ghat{{\widehat{G}}}
\def\Hhat{{\widehat{H}}}

\def\Bra{{\brr_{\rm a}\,}}
\def\Bro{{\brr_1\,}}

\def\Xbar{{\overline{X}}}
\def\UPicc{{\rm UPic}}
\def\UPic{{\rm KD}'}
\def\Gal{{\rm Gal}}
\def\Picc{{\rm Pic}}
\def\kbar{{\bar{k}}}

\def\sH{{\mathcal{H}}}

\def\Phat{{\widehat{P}}}

\def\sO{{\mathcal{O}}}
\def\sK{{\kbar}}

\def\OXbar{{\kbar[\Xbar]^*}}
\def\KXbar{{\kbar(\Xbar)^*}}

\def\Zalg{Z_\alg}
\def\Div{{\rm Div}}
\def\divisor{{\rm div}}
\def\id{{\rm id}}
\def\into{{\hookrightarrow}}

Dans cet  appendice nous prouvons  la proposition \ref{prop:main-prop} que nous utilisons  dans la preuve du lemme \ref{lem:Z}.

Soit $T$ un tore sur un corps $k$.
On pose $G=T$,  $H=1$, $X=T$.
Alors \cite[thm.~7.2]{borvh2} donne un isomorphisme
$$
\muB\colon \Bra T\isoto H^1(k,[\Ghat\to \Hhat])= H^1(k,[\That\to 0]=H^2(k,\That),
$$
o\`u $\Bra T=\Bro T/\Br k$.

\begin{prop}\label{prop:main-prop}
Le diagramme suivant est anti-commutatif:
$$
\xymatrix{
T(k)\ar@{=}[d] &\times &\Bra T\ar[d]^\muB\ar[r]^-{\langle,\rangle}   &\Br k\ar@{=}[d]\\
T(k)           &\times  &H^2(k,\That)\ar[r]^-\cup   &\Br k
}
$$
o\`u l'accouplement sup\'erieur est d\'efini  en utilisant les homomorphismes
$$
\Bra T\isoto\Brr_{1,e}\, T\into \Bro T.
$$
\end{prop}

\def\Gm{{\mathbf{G}_m}}
\def\trG{{\R p_* {\bf G}_{m,X}}}

Pour prouver la proposition nous avons besoin d'une s\'erie de lemmes.

\begin{lem}\label{lem:abc}
Soit $X$ une $k$-vari\'et\'e lisse et g\'eom\'etriquement int\`egre
munie d'un $k$-point $m$. On note $\kbar[X]^*:=H^0(\ov X,\G)$
le groupe multiplicatif des fonctions inversibles sur $X$ et
$U(\ov X):=\kbar[X]^*/\kbar^*$.
Alors le diagramme suivant est commutatif:
\begin{equation*}
\xymatrix{
H^2(k,U(\Xbar))\ar[r]^\beta\ar[rd]_\gamma &\Bra X\ar[d]^\alpha  \\
                           &H^1(k,\UPic(X)),
}
\end{equation*}
o\`u $\beta$ est l'homomorphisme canonique de  Sansuc \cite[Lemme 6.3(ii)]{sansuc},
$\alpha$ est l'isomorphisme canonique  de \cite[Cor. 2.20(i)]{borvh},
et $\gamma$ est l'homomorphisme canonique induit  par le morphisme  de complexes
$$
\deltaN\colon U(\Xbar)=\sH^{0}(\UPic(X)[-1])\into\UPic(X)[-1].
$$
\end{lem}

\dem
(1) Rappelons 
que $\trG$ est un certain complexe $[C_0\to C^1][-1]$
de $\Gal(\kbar/k)$-modules en degr\'es 0 et 1
(d\'efini \`a un isomorphisme canonique pr\`es dans la cat\'egorie deriv\'ee)
avec 0-cohomologie $\sH^0(\trG)=\kbar[X]^*$.
Le groupe $\kbar^*$ s'injecte canoniquement  dans $\kbar[X]^*$,
et on a
$$
\UPic(X):=[C^0/\kbar^*\to C^1].
$$
La cohomologie en degree $-1$ du complexe $\UPic(X)$ est donn\'ee par
$$
\sH^{-1}(\UPic(X))=U(\Xbar):=\kbar[X]^*/\kbar^*.
$$
On note qu'il y a un isomorphisme canonique
$$\alpha'\colon\Bro X\isoto H^2(k,\trG),$$
voir \cite[2.18]{borvh}.

(2) On consid\`ere  l'homomorphisme
$$
H^2(k,\kbar[X]^*)\to H^2(k,\trG)\labelto{(\alpha')^{-1}}\Bro X,
$$
o\`u le premier homomorphisme est induit par le morphisme de complexes
$$
\kbar[X]^* \into \trG,
$$
alors il r\'esulte des d\'efinitions que la composition co\"\i ncide avec l'homomorphisme
de Sansuc \cite[lemme 6.3(i)]{sansuc}.

(3) On consid\`ere l'homomorphisme
$$
\beta\colon H^2(k,U(\Xbar))\to\Bra X
$$
de Sansuc \cite[lemme 6.3(ii)]{sansuc}.
Cet homomorphisme est d\'efini comme suit:
le $k$-point $m\in X(k)$ d\'efinit une  section homomorphique  $\Gal(\kbar/k)$-\'equivariante
$\sigma_{m}\colon U(\Xbar)\to \kbar[X]^*$ de l'\'epimorphisme canonique
$\kbar[X]^*\to U(\Xbar)$ donn\'ee par
$$
\sigma_{m} ([f])(x)=f(x)/f(m)
$$
pour la classe $[f]$ d'une fonction $f\in \kbar[X]^*$.
Alors $\beta$ est la composition
$$
H^2(k,U(\Xbar))\labelto{(\sigma_{m})_*} H^2(k,\kbar[X])\to H^2(k, \trG)\labelto{(\alpha')^{-1}}
\Bro X\to \Bra X\, .
$$

(4) On a un  diagramme commutatif
$$
\xymatrix{
\kbar[X]^*\ar[d]\ar[r]    &\trG \ar[d] \\
U(\Xbar)\ar@{.>}@/_1pc/[u]_{\sigma_{m}}  \ar[r]  &\UPic(X)[-1],
}
$$
qui  induit un diagramme commutatif d'hypercohomologie
$$
\xymatrix{
H^2(k, \kbar[X]^*)\ar[d]_\rho\ar[r]^-{\gamma'}       &H^2(k,\trG)\ar[d]       &\Bro X\ar[d]^-{\rho'}\ar[l]_-{\alpha'} \\
H^2(k,U(X))\ar@{.>}@/_1pc/[u]_{(\sigma_{m})_*}\ar[r]^-\gamma   &H^1(k,\UPic(X))   &\Bra X,\ar[l]_-\alpha
}
$$
Comme $\beta=\rho'\circ(\alpha')^{-1}\circ\gamma'\circ(\sigma_{m})_*$\,,
on en conclut que $\gamma=\alpha\circ\beta$.
\enddem

\begin{lem}\label{lem:minus}
Soient $T$, $G$, $H$ et $X$ comme au d\'ebut de l'appendice;
on d\'esigne par
$$
\deltaB\colon\That\to [\KXbar/\kbar^*\to\Div(\Xbar)][-1]
$$
l'isomorphisme dans le cat\'egorie d\'eriv\'ee de \cite[Thm.~5.8]{borvh2}.
Alors $\deltaB$ est l'oppos\'ee du quasi-isomorphisme naturel
$$
\delta_{\rm nat}\colon\That\to [\KXbar/\kbar^*\to\Div(\Xbar)][-1]
$$
induit par l'injection naturelle $\That=U(\Xbar)\into \KXbar/\kbar^*$,
\`a savoir $\deltaB=-\delta_{\rm nat}$.
\end{lem}

\dem
On commence par consid\`erer \cite[diagramme (13) apr\`es le Thm.~4.10]{borvh2}:
\begin{equation*}\label{eq:diagram-isomorphism}
 \xymatrix{
{\Zalg^1(\Gbar, \OXbar)}\ar[d]^\sigma  &{\Zalg^1(\Gbar, \OXbar)}\oplus\sK(\Xbar)^*/\kbar^*\ar[d]^\psi\ar[l]\ar[r]
                                                            &\sK(\Xbar)^*/\kbar^*\ar[d]^{\divisor}\\
\Picc_G\,\Xbar                          &\UPicc_G(\Xbar)^1\ar[l]\ar[r]                  &\Div(\Xbar)
}
\end{equation*}
(voir \cite{borvh2} pour les notations),
o\`u la fl\`eche $\psi$ est donn\'ee par
$$
\psi(z,[f])=(z\cdot d^0(f),\divisor(f))\in \UPicc_G(\Xbar)^1\subset  \Zalg^1(\Gbar,\KXbar)\oplus\Div(\Xbar),
$$
et toutes les fl\`eches non \'etiquet\'ees sont les fl\`eches \'evidentes.
\medskip

On a $\Zalg^1(\Gbar,\OXbar)=\That$, $\Picc_G\,\Xbar=\Hhat=0$.
On pose
$$
a=\id\colon \That\isoto \Zalg^1(\Gbar,\OXbar).
$$
Clairement, $a$ induit un morphisme de complexes
$$
a'\colon \That\to [\Zalg^1(\Gbar,\OXbar)\to\Picc_G\,\Xbar][-1].
$$
On consid\`ere l'injection canonique
$$
c\colon \That=U(\Xbar)\into \KXbar/\kbar^*.
$$
Clairement, $c$ induit un morphisme de complexes
$$
c'=\deltaN\colon \That \to [\KXbar/\kbar^*\to\Div(\Xbar)][-1].
$$
On d\'efinit
$$
b\colon \That \to \Zalg^1(\Gbar,\OXbar) \oplus  \KXbar/\kbar^*
$$
par $b(\chi)=(a(\chi), -c(\chi))$ pour $\chi\in\That$.
Alors il r\'esulte de la  formule pour $\psi$ que $\psi\circ\beta=0$
(on utilise la formule pour $b$ avec le signe oppos\'e \`a $c(\chi)$!),
donc $b$ induit un morphisme de complexes
$$
b'\colon \That\to [\Zalg^1(\Gbar,\OXbar) \oplus  \KXbar/\kbar^* \to \UPicc_G(\Xbar)^1][-1].
$$
Clairement, on a un diagramme commutatif
\begin{equation*}\label{eq:diagram-That}
 \xymatrix{
\That\ar[d]^a\ar[r]^{=}     &\That\ar[d]^b\ar[r]^{=}    &\That\ar[d]^{-c}  \\
{\Zalg^1(\Gbar, \OXbar)}\ar[d]^\sigma  &{\Zalg^1(\Gbar, \OXbar)}\oplus\sK(\Xbar)^*/\kbar^*\ar[d]^\psi\ar[l]\ar[r]
                                                            &\sK(\Xbar)^*/\kbar^*\ar[d]^{\divisor}\\
\Picc_G\,\Xbar                          &\UPicc_G(\Xbar)^1\ar[l]\ar[r]                  &\Div(\Xbar).
}
\end{equation*}

La partie inf\'erieure du diagramme induit des isomorphismes
(dans la cat\'egorie d\'eriv\'ee)~:
$$[\Zalg^1(\Gbar, \OXbar) \to \Picc_G\,\Xbar] \simeq
[\Zalg^1(\Gbar, \OXbar)\oplus\sK(\Xbar)^*/\kbar^* \to \UPicc_G(\Xbar)^1]
$$
et
$$[\Zalg^1(\Gbar, \OXbar)\oplus\sK(\Xbar)^*/\kbar^* \to \UPicc_G(\Xbar)^1]
\simeq [\kbar(X)^*/\kbar^* \to \div \ov X]]$$
Or l'isomorphisme (dans la cat\'egorie d\'eriv\'ee) $\deltaB$
est la composition de $a'$ et du compos\'e des deux isomorphismes ci-dessus.
On d\'eduit du diagramme que $\deltaB=-c'$.
Comme  $\deltaN=c'$, on conclut  que $\deltaB=-\deltaN$.
\enddem

\begin{lem}\label{lem:main-lemma}
On d\'esigne par $\lambdaS\colon H^2(k,\That) \isoto\Bra T$ l'isomorphisme canonique de Sansuc \cite[lemme 6.9]{sansuc}.
Alors
$\muB\circ\lambdaS=-\id\colon H^2(k,\That)\to \Bra T\to  H^2(k,\That)$.
\end{lem}

\dem
On a un  diagramme commutatif
$$
\xymatrix{
H^2(k,\That)\ar@{.>}[d] \ar[r]^{\beta=\lambda_S}  &\Bra X\ar[d]^\alpha \ar[ld]_\muB  \\
H^2(k,\That)\ar[r]^-{(\deltaB)_*}              &H^1(k,\UPic(X))
}
$$
Par d\'efinition $\lambda_S=\beta$, $\muB=(\deltaB)_*^{-1}\circ\alpha$.
D'apr\`es le lemme \ref{lem:minus}, on a $\deltaB=-\deltaN$.
Par d\'efinition $\gamma=(\deltaN)_*$, voir le lemme \ref{lem:abc}.
Donc
$$
\muB=(-\deltaN)_*^{-1}\circ\alpha=-(\deltaN)_*^{-1}\circ\alpha=-\gamma^{-1}\circ\alpha.
$$
Ainsi
$$
\muB\circ\lambdaS=-\gamma^{-1}\circ\alpha\circ\beta=-\id,
$$
parce que d'apr\`es le lemme \ref{lem:abc}, on a $\gamma=\alpha\circ\beta$.
\enddem

 \paragraph{Preuve de la proposition \ref{prop:main-prop} :}
On consid\`ere le diagramme
\begin{equation*}\label{eq:big-diag}
\xymatrix{
T(k)\ar@{=}[d] &\times &H^2(k,\That)\ar[d]^\lambdaS\ar[r]^-\cup   &\Br k\ar@{=}[d]\\
T(k)\ar@{=}[d] &\times &\Bra T\ar[d]^\muB\ar[r]^-{\langle,\rangle}   &\Br k\ar@{=}[d]\\
T(k)           &\times  &H^2(k,\That)\ar[r]^-\cup                  &\Br k
}
\end{equation*}
Le rectangle sup\'erieur de ce diagramme est commutatif, gr\^ace \`a
\cite[diagram (8.11.2)]{sansuc}.
D'autre part,  si on supprime  la ligne m\'ediane de ce diagramme,
alors le diagramme qui en r\'esulte sera
anti-commutatif, parce que d'apr\`es le lemme \ref{lem:main-lemma} $\muB\circ\lambdaS=-\id$.
On voit que le  rectangle inf\'erieur du  diagramme est anti-commutatif.
\enddem

\section{Une autre preuve du th\'eor\`eme~\ref{compatible}}\label{AppDem}

Dans cet appendice, nous donnons une preuve de l'expression
explicite de $KD'(X)$ et de la formule de
compatiblit\'e (th\'eor\`eme~\ref{compatible}) ind\'ependante
de \cite{borvh2}, sous l'hypoth\`ese que le stabilisateur $H$ est de
type (ssumult). Au passage, on \'etablit quelques r\'esultats
g\'en\'eraux sur le complexe $KD'$ associ\'e \`a un torseur sous un
groupe de type multiplicatif.

\smallskip Dans toute la suite, $k$ est un corps et
$p : X \to \spec k$ d\'esigne 
une vari\'et\'e lisse et g\'eom\'etriquement int\`egre
sur $k$.

\begin{lem} \label{dualcomplex}
Soit $[S \to T]$ un complexe de $k$-groupes de type multiplicatif. Soit
$[\widehat T \to \widehat S]$ le complexe dual. Alors on a
$$H^0_{\fp}(X,[S \to T])=\Hom_k([\widehat T \to \widehat S], KD(X))$$
o\`u $\Hom_k(...)$ d\'esigne les homomorphismes dans la cat\'egorie
d\'eriv\'ee born\'ee ${\cal D}(k)$ des faisceaux \'etales sur $\spec k$.
\end{lem}

Pour simplifier les notations, on \'ecrira souvent $H^0(X,[S \to T])$ pour
$H^0(X,[p^*S \to p^* T])$.

\dem L'argument est le m\^eme que
dans la proposition 1.1 de \cite{opendesc}. Pour tout complexe born\'e
$C^{\bullet}$ (resp. $D^{\bullet}$) de faisceaux \'etales sur $\spec k$
(resp. sur $X$), on a
$\Hom_k(C^{\bullet},p_* D^{\bullet})=\Hom_X(p^* C^{\bullet},
D^{\bullet})$.
On en d\'eduit
que le foncteur $\R {\rm Hom}_X([p^* \widehat T \to p^* \widehat S],.)$
est le compos\'e
des foncteurs $\R p_*$ et $\R \Hom_k([\widehat T \to \widehat S])$,
ce qui permet d'obtenir (pour tout $n \geq 0$)~:
$$\R ^n \Hom_k([\widehat T \to \widehat S],\R  p_* {\bf G}_{m,X})=
\R^n \Hom_X([p^* \widehat T \to p^* \widehat S],{\bf G}_{m,X}) \, .$$
Comme par d\'efinition
$KD(X)=\tau_{\leq 1} \R p_* {\bf G}_{m,X}[1]$, le triangle exact
$$\tau_{\leq 1} \R p_* {\bf G}_{m,X} \to \R p_* {\bf G}_{m,X} \to
\tau_{\geq 2} \R p_* {\bf G}_{m,X} \to \tau_{\leq 1} \R p_* {\bf G}_{m,X}[1]$$
donne d'autre part
$$\Hom_k([\widehat T \to \widehat S], KD(X))=
\R ^1 \Hom_k([\widehat T \to \widehat S],\R p_* {\bf G}_{m,X})$$
via le fait que $\tau_{\geq 2} \R p_* {\bf G}_{m,X}$
est acyclique en degr\'e $<2$.
Pour conclure il suffit alors de montrer la formule
$$ \R^1 \Hom_X([p^* \widehat T \to p^* \widehat S],{\bf G}_{m,X})=
H^0_{\fp}(X,[S \to T])$$
qu'on peut r\'ecrire
$$H^0_{\fp}(X,[S \to T])=
\Hom_X([p^* \widehat T \to p^* \widehat S], {\bf G}_{m,X}[1])$$
ce qui r\'esulte du lemme~\ref{defpairing}.

\enddem

\rem L\`a encore, la m\^eme m\'ethode montre que l'assertion reste
vraie si on remplace
$\Hom_k(...)$ par $\Hom_{{\mathcal D}_{\fp}(k)}(...)$ dans le groupe
de droite.

\begin{prop} \label{genprop}
Soit $[S \to T]$ un complexe de $k$-groupes de type multiplicatif.
Soit $f : Y \to X$ un $S$-torseur. On suppose de plus que $Y$ est muni
d'une trivialisation du torseur $Y \wedge^S T \to X$ obtenu par
changement de groupe structural. Soit $[[Y]]$ la classe correspondante
dans $H^0_{\fp}(X,[S \to T])$ et
$u \in \Hom_k([\widehat T \to \widehat S], KD(X))$
le morphisme associ\'e comme dans le lemme~\ref{dualcomplex}.
On d\'efinit un morphisme $\mu : [\widehat T \to \widehat S] \to
KD'(X)$ dans ${\cal D}(k)$ en composant $u$ avec le morphisme
canonique $v : KD(X) \to KD'(X)$. Alors on a,
pour tout $a \in H^1(k,[\widehat T \to \widehat S])$,
la formule $$r(p^* a \cup [[Y]])=\mu_*(a) \, .$$
\end{prop}

\dem C'est tout \`a fait similaire \`a la preuve du th\'eor\`eme 1.4.
de \cite{opendesc}. On a (via le lemme~\ref{dualcomplex} et sa preuve)
un diagramme commutatif
{\small
$$\begin{array}{ccccc}
H^1_{\fp}(k,[\widehat T \to \widehat S])& \times& H^0_{\fp}(X,[S \to T])&\to&\br(X)\\
||&&||&& \uparrow \\
H^1_{\fp}(k,[\widehat T \to \widehat S])& \times&
\Hom_{{\mathcal D}_{\fp}(k)}([\widehat T \to \widehat S],
\R p_* {\bf G}_{m,X}[1])&\to& H^1_{\fp}(k,\R p_* {\bf G}_{m,X}[1])\\
||&&||&&\uparrow\\
H^1(k,[\widehat T \to \widehat S])& \times&
\Hom_k([\widehat T \to \widehat S], KD(X))&\to&H^1(k,KD(X))
\end{array}$$
}
o\`u l'accouplement de la premi\`ere ligne est donn\'e par le cup-produit
(cf. \cite{milne}, Prop. V.1.20). On obtient donc,
une fois qu'on a identifi\'e $\bra X$ avec $H^1(k,KD(X))$~:
$$p^* a \cup [[Y]]=u_*(a) \, .$$
Ainsi
$$r(p^* a \cup [[Y]])=r(u_*(a))=v_*(u_*(a))=\mu_*(a)$$
par d\'efinition de $\mu$ et de $r$.

\enddem

Rappelons la d\'efinition suivante (\cite{opendesc}, section 2)~:

\begin{defi}
{\rm
Soit $V$ une $k$-vari\'et\'e (pas forc\'ement int\`egre) munie d'une action
$m : S \times V \to V$ d'un $k$-groupe de type multiplicatif $S$. On d\'efinit
$\kbar [V]^*_S$ comme le sous-groupe de $\kbar [V]^*:=H^0(\ov V, \G)$
constitu\'e des fonctions $f$ pour lesquelles il existe un caract\`ere
$\chi \in \widehat S$ v\'erifiant $m^* f=\chi.f$ dans
$H^0(\ov S \times \ov V,\G)$. De m\^eme si $\pi : Y \to X$ est un $S$-torseur,
on d\'efinit le faisceau \'etale $(\pi_* {\bf G}_{m,Y})_S$ comme le sous-faisceau
de $\pi_* {\bf G}_{m,Y}$ tel que pour tout morphisme \'etale $U \to X$,
le groupe $(\pi_* {\bf G}_{m,Y})_S(U)$ consiste en les $f : Y_U \to {\bf G}_{m,U}$
telles qu'il existe un caract\`ere $\chi :
S_U \to {\bf G}_{m,U}$ avec $m^* f=\chi.f$,
o\`u $Y_U:=Y \times_k U$ et $S_U:=S \times_k U$.
}
\end{defi}

Noter que si
l'action du sch\'ema en groupes $S$ sur $U(\ov V)=\kbar[V]^*/\kbar^*$
est triviale (i.e. la fl\`eche $m^* : H^0(\ov V,\G)/\kbar^* \to
H^0(\ov S \times \ov V,\G)/\kbar^*$ induite par l'action $m : S \times V
\to V$ co\"{\i}ncide avec la fl\`eche induite par la projection
$\ov S \times \ov V \to \ov V$), alors $\kbar [V]^*_S=\kbar [V]^*$.
C'est en particulier le cas si $S$ est un tore et
$V$ est g\'eom\'etriquement connexe (\cite{opendesc}, remarque page 9).

\begin{prop}
\label{prop suite 5 termes}
Soit $S$ un $k$-groupe de type multiplicatif lisse. Soit $\pi : Y \to X$ un
$S$-torseur avec $Y$ g\'eom\'etriquement connexe.
On suppose que l'action du sch\'ema en groupes $S$ sur
$U(\ov Y)=\kbar[Y]^*/\kbar^*$ est triviale.

Alors on a une suite exacte naturelle
$$0 \to U(\ov X) \to U(\ov Y) \to \widehat S \to \Pic(\ov X) \to \Pic(\ov Y) \, .$$
o\`u la fl\`eche $ \widehat S \to \Pic(\ov X)$ est le type du torseur $\pi$.
\end{prop}

\dem On peut supposer $k=\kbar$. L'hypoth\`ese que $S$ (et donc aussi $F$) est
lisse implique alors que le $k$-groupe fini (et commutatif) $F$ est
constant.
On note $T$ la composante connexe de l'identit\'e dans $S$ et
on pose $F := S / T$ et $Z := Y/T$.
On a alors un diagramme commutatif de torseurs :
\begin{displaymath}
\xymatrix{
Y \ar[rd]^T \ar[dd]^S & \\
& Z \ar[ld]^F \\
X & \, .
}
\end{displaymath}
On en d\'eduit le diagramme commutatif suivant :
\begin{displaymath}
\xymatrix{
0 \ar[d] & & \\
\widehat F \ar[d] \ar[r]^= & \widehat F \ar[d] &  \\
\widehat S \ar[d] \ar[r] & \Pic X \ar[r] \ar[d] & \Pic Y \ar[d]^= \\
\widehat T \ar[r] \ar[d] & \Pic Z \ar[r] & \Pic Y \\
0 & & \, .
}
\end{displaymath}
La premi\`ere colonne est exacte, la ligne inf\'erieure est exacte par \cite{sansuc}, proposition 6.10.
L'exactitude de la deuxi\`eme colonne
est donn\'ee par la suite exacte des termes de
bas degr\'e associ\'ee \`a la suite spectrale de
Hochschild-Serre~: $E_2^{p,q} := H^p(F, H^q(Z, {\bf G}_{m}))
\Rightarrow H^{p+q}(X, {\bf G}_{m})$
du rev\^etement galoisien $Z \to X$.
En effet
$H^0(F,H^1(Z, {\bf G}_{m}))$ est un sous-groupe de $\pic Z$, et
$H^1(F,H^0(Z, {\bf G}_{m}))$ est un quotient de $H^1(F,k^{*})=\widehat F$
via la suite exacte
$$0 \to k^* \to H^0(Z, {\bf G}_{m}) \to U(Z) \to 0 \, ,$$
et l'hypoth\`ese que le groupe fini $F$ agit trivialement sur
$U(Y)$, donc aussi sur $U(Z)$ (qui est un
groupe ab\'elien libre de type fini car $Z$ est lisse et connexe),
ce qui implique $H^1(F,U(Z))=0$.

\smallskip

Une chasse au diagramme utilisant l'exactitude des deux premi\`eres colonnes et de la derni\`ere ligne
assure alors que la ligne centrale est exacte.
Cela conclut la preuve, gr\^ace \`a la proposition 2.5, (i) et (iii), de \cite{opendesc}.
\enddem

On rappelle le lemme suivant d'alg\`ebre homologique (voir par exemple \cite{DemBrauer}, lemme 2.3) :

\begin{lem}
\label{lem homologique}
Soient $\cal A$, $\cal B$ deux cat\'egories ab\'eliennes et $F : \cal A \to \cal B$ un foncteur exact \`a gauche.
On suppose que $\cal A$ a suffisamment d'injectifs.
Soit $0 \to A \to B \to C \to 0$ une suite exacte d'objets de $\cal A$. Si le cobord $F(C) \to R^1 F(A)$ est un \'epimorphisme,
alors on a un isomorphisme canonique dans la cat\'egorie d\'eriv\'ee :
$$[F(B) \to F(C)] \to \tau_{\leq 1}(\R F)(A) \, .$$
\end{lem}

\dem
C'est le lemme 2.3 de \cite{DemBrauer}, dans le cas $n=1$.
\enddem

\begin{lem} \label{torstore}
Soit $S$ un $k$-groupe de type multiplicatif lisse. Soit $\pi : Y \to X$ un
$S$-torseur, avec $Y$ g\'eom\'etriquement connexe.
On note $p : X \to \spec k$ le morphisme structural
et $q=(p \circ \pi) : Y \to \spec k$.
On fait les deux hypoth\`eses suivantes~:

\smallskip

(a) L'action du sch\'ema en groupes $S$ sur $U(\ov Y)=\kbar[Y]^*/\kbar^*$
est triviale.

\smallskip

(b) On a $\pic \ov Y=0$.

\smallskip

Alors on a un isomorphisme dans ${\cal D}(k)$~:
\begin{equation} \label{trikd}
[\kbar[Y]^* \to \widehat S] \to KD(X)
\end{equation}
o\`u la fl\`eche induite $$\widehat S \to KD(X)$$ s'identifie \`a la classe
$[Y] \in H^1_{\fp}(X,S)= \Hom_k(\widehat S,KD(X))$ (au signe pr\`es).
\end{lem}

\dem On a une suite exacte de faisceaux \'etales sur $X$~:
\begin{equation} \label{rosensuite}
0 \to {\bf G}_{m,X} \to (\pi_* {\bf G}_{m,Y})_S \to p^* \widehat S \to 0
\end{equation}
(voir \cite{opendesc}, proposition 2.5(i)). Notons que $p_* p^* \widehat S=\widehat S$ par connexit\'e de $\ov X$.
On applique le foncteur $p_ *$ : la suite exacte de la proposition \ref{prop suite 5 termes} et l'hypoth\`ese (b) assurent
que l'on a une suite exacte de modules galoisiens
$$0 \to p_* {\bf G}_{m,X} \to p_* (\pi_* {\bf G}_{m,Y})_S \to \widehat S \to R^1 p_* {\bf G}_{m,X} \to 0 \, .$$
Alors le lemme \ref{lem homologique} assure que le morphisme naturel
$$[p_* (\pi_* {\bf G}_{m,Y})_S \to \widehat S] \to \tau_{\leq 1} (\R p_*) {\bf G}_{m,X} = KD(X)$$
est un isomorphisme. Or l'hypoth\`ese (a) assure que $p_* (\pi_* {\bf G}_{m,Y})_S = \kbar[Y]^*$, d'o\`u finalement un isomorphisme
$$[\kbar[Y]^* \to \widehat S] \to KD(X) \, .$$

\smallskip

Comme la classe $[Y] \in
\ext^1_X(\widehat S, {\bf G}_{m,X})=\Hom_X(\widehat S, {\bf G}_{m,X}[1])$
correspond au signe pr\`es (\cite{opendesc}, Prop. 2.5(ii))
\`a l'extension donn\'ee
par (\ref{rosensuite}), la fl\`eche induite $\alpha : \widehat S \to KD(X)$
correspond bien \`a la classe $[Y]$ (cf. \cite{opendesc}, appendice B).

\enddem

\begin{prop} \label{torsmult}
On garde les hypoth\`eses et notations du lemme~\ref{torstore}.
Soit $T$ le $k$-tore dont le module des caract\`eres est
$\widehat T=\kbar[Y]^* /\kbar^*$ (qui est muni d'un
morphisme canonique $\widehat T \to \widehat S$ par
\cite{opendesc}, Prop. 2.5).

\smallskip

(a) On a un
isomorphisme dans ${\cal D}(k)$~:
\begin{equation} \label{triafini}
[\widehat T \to \widehat S] \to KD'(X) \, .
\end{equation}

\smallskip

(b) On suppose de plus que $Y$ poss\`ede un point rationnel
$y \in Y(k)$.
Alors $y$ induit une trivialisation de $Y \wedge^S T$, telle que
le morphisme $\mu$ associ\'e \`a $Y$ comme dans la proposition~\ref{genprop}
soit un isomorphisme dans ${\cal D}(k)$.
\end{prop}

\dem (a)
L'assertion r\'esulte du lemme~\ref{torstore} : on a un isomorphisme $\mu_0 : [\kbar[Y]^* \to \widehat{S}] \to KD(X)$
qui induit un isomorphisme $\overline{\mu}_0 : [\widehat T \to \widehat S] \to KD'(X)$.

\smallskip

(b) On dispose du diagramme commutatif (fl\`eches rectilignes) suivant :
\begin{displaymath}
\xymatrix{
\textup{Hom}_k([\kbar[Y]^* \to \widehat S], KD(X)) \ar[r] \ar@/^2pc/[d]^{s_y} & Hom_k(\widehat S, KD(X)) \\
\textup{Hom}_k([\widehat T \to \widehat S], KD(X)) \ar[u] \ar[r] \ar[d]^{\cong} & Hom_k(\widehat S, KD(X)) \ar[d]^{\cong} \ar[u]^= \\
H^0(X, [S \to T]) \ar[r] & H^1(X,S) \, ,
}
\end{displaymath}
o\`u $s_y$ est induite par la section de $\kbar[Y]^* \to \kbar[Y]^*/\kbar^*$ donn\'ee par $y \in Y(k)$.

On note $\mu_{0,y} := s_y(\mu_0) : [\widehat T \to \widehat S] \to KD(X)$.

Le diagramme pr\'ec\'edent et la seconde partie du lemme \ref{torstore} assurent que l'image de $\mu_{0,y}$ dans $H^0(X,[S \to T])$
est la classe du torseur $Y \to X$ muni d'une trivialisation du torseur $Y \wedge^S T$,
de sorte que le morphisme $[\widehat T \to \widehat S] \to KD'(X)$ associ\'e \`a cette trivialisation
(comme \`a la proposition \ref{genprop}) soit exactement $\overline{\mu}_0$.
Cela assure le point (b) de la proposition.

\enddem

\begin{prop} \label{homog}
Soit $G$ un $k$-groupe lin\'eaire, r\'eductif et connexe.
Soit $H$ un $k$-sous-groupe de type (ssumult) de $G$, on pose
$S=H^{\mult}$ et $H'=\ker [H \to H^{\mult}]$.
On suppose de plus que $\pic \ov G=0$  et on pose $T=G^{\tor}$.
On pose $X=G/H$ et
on note  $[[G]] \in H^0_{\fp}(X,[H \to G])$ la
classe induite par la trivialisation
du $X$-torseur (sous $G$) $G \wedge^H G$ donn\'ee par le neutre $e \in G(k)$.
On note \'egalement $$\tau : X(k) \to H^0_{\fp}(k,[H \to G])
\quad x \mapsto [[G]](x)$$
la bijection donn\'ee par l'\'evaluation.

\smallskip

Soit $Y:=G/H'$, qu'on peut voir comme un $X$-torseur
sous $S$ muni d'une trivialisation de $Y \wedge^S T$.
Alors la fl\`eche $$\mu : [\widehat T \to \widehat S] \to
KD'(X)$$ associ\'ee (comme dans la proposition~\ref{genprop})
\`a $Y$ est un isomorphisme dans ${\cal D}(k)$.

\end{prop}

\dem On observe
que comme $$\kbar[Y]^*/\kbar^*=\kbar[G]^*/\kbar^*=\widehat G=\widehat T$$
par le lemme de Rosenlicht, l'action de $S$ sur $\kbar[Y]^*/\kbar^*$ est
triviale. Comme $H$ est de type (ssumult), le groupe $H'$ est connexe,
lisse, et sans caract\`eres; d'autre part $\ov G$ est une vari\'et\'e
$\kbar$-rationnelle avec $\pic \ov G=0$.
Par \cite{sansuc}, Prop. 6.10, on a alors $\pic \ov Y=0$.
Il suffit alors
d'appliquer la proposition~\ref{torsmult}(b) au $S$-torseur $Y \to X$.

\enddem

\rem Le m\^eme argument marche avec un $H$ un peu plus g\'eneral
(extension d'un groupe de type multiplicatif lisse par un groupe
$H'$ lisse, de composante neutre $H'_0$  v\'erifiant~: $H'/H'_0$ et
$H'_0$ sont sans caract\`eres).

\begin{theo}
On garde les hypoth\`eses et notations de la proposition~\ref{homog}.
Soit $$\ab : H^0_{\fp}(k,[H \to G]) \to H^0_{\fp}(k,[S \to T])$$ l'application
d'ab\'elianisation et $\ab^0:=(\ab \circ \tau) : X(k) \to
H^0_{\fp}(k,[S \to T])$.

Soit $\Phi_X : \brr_{1,e} X \to H^1(k,[\widehat T \to
\widehat S])$ l'isomorphisme
d\'efini par $\Phi_X=(\mu_*)_{-1} \circ r$. Alors on a (au signe
pr\`es)~:
$$\phi_X(\alpha) \cup \ab^0(x)=\alpha(x)$$
pour tout $x \in X(k)$ et tout $\alpha \in \brr_{1,e} X$.
\end{theo}

\dem D'apr\`es la proposition~\ref{genprop}, on a
$$\Phi_X(p^* a \cup [[Y]])=\mu_* ^{-1}(r(p^* a \cup [[Y]]))=a$$
pour tout $a \in H^1(k,[\widehat T \to \widehat S])$.
Posons $a=\Phi_X(\alpha)$. Comme $\Phi_X$ est injective, on obtient
$$p^*(\Phi_X(\alpha)) \cup [[Y]]=\alpha$$
d'o\`u en \'evaluant en $x$ et en remarquant que
$[[Y]](x)=\ab(\tau(x))= \ab^0(x)$~:
$$\Phi_X(\alpha) \cup \ab^0(x)=\alpha(x)$$
qui est la formule voulue.

\enddem

\merci Les auteurs tiennent \`a remercier chaleureusement 
J-L. Colliot-Th\'el\`ene, C.D. Gonz\'alez-Avil\'es
et B. Kahn pour d'int\'eressantes discussions pendant la pr\'eparation de
cet article.

\bigskip

\noindent
Raymond and Beverly Sackler  School of Mathematical Sciences, Tel Aviv University,
69978 Tel Aviv, Israel
\smallskip

\noindent
borovoi@post.tau.ac.il

\bigskip 

\noindent
Institut de Math\'ematiques de Jussieu (IMJ), Universit\'e Pierre et Marie Curie, 4 place Jussieu, 
75252 Paris Cedex 05, France
\smallskip

\noindent
demarche@math.jussieu.fr

\bigskip 

\noindent 
Universit\'e Paris-Sud, Laboratoire de math\'ematiques, B\^atiment 425,
F-91405 Orsay Cedex, France 

\smallskip

\noindent
david.harari@math.u-psud.fr

\end{document}